\documentclass{article}
\usepackage{amssymb}
\usepackage{amsfonts}
\usepackage{amsmath}
\usepackage[left=1.5in, right=1.5in, bottom=1.5in,top=1.5in]{geometry}

\input{tcilatex}
\begin{document}

\author{ M. El- Morshedy\thanks{%
Corresponding author. Permanent address: Department of Mathematics, College
of Science, Mansoura University, Mansoura 35516, Egypt. E-mail address:
Mah\_elmorshedy@mans.edu.eg (M. El-Morshedy).} and A. A. Khalil \\
Mathematics Department, \\
Faculty of Science, Mansoura University, Mansoura,\\
Egypt.}
\title{\textbf{Bivariate Discrete Exponentiated Weibull Distribution:
Properties and Applications}}
\date{}
\maketitle

\begin{abstract}
In this paper, a new bivariate discrete distribution is introduced which
called bivariate discrete exponentiated Weibull (BDEW) distribution. Several
of its mathematical statistical properties are derived such as the joint
cumulative distribution function, the joint joint hazard rate function,
probability mass function, joint moment generating function, mathematical
expectation and reliability function for stress-strength model. Further, the
parameters of the BDEW distribution are estimated by the maximum likelihood
method. Two real data sets are analyzed, and it was found that the BDEW
distribution provides better fit than other discrete distributions.

\end{abstract}

\ \ \ 

\section{Introduction}

The Weibull distribution (1951) is one of the most important and
well-recognized continuous probability model in research and also in
teaching. It has become important because it can be able to assume the
properties of many varies types of continuous distributions. And so, if you
have right-skewed data, left-skewed data or symmetric data, you can use
Weibull to model it. Moreover, the hazard rate of its can be constant,
increasing or decreasing. That flexibility of Weibull distribution has made
many researchers using it into their data analysis in different fields such
as medicine, pharmacy, Engineering, astronomy, electronics, reliability,
industry, space science, social sciences, economics and environmental. In
previous years, many researchers interested in the distributions theory have
provided many generalizations or extensions of the Weibull distribution.
See, Mudholkar and Srivastara (1993), Bebbington et al. (2007), Sarhan and
Apaloo (2013), El-Gohary et al. (2015), El- Bassiouny et al. (2017), El-
Morshedy et al. (2017), among others.

Despite the great importance for the continuous probability distributions,
there are many practical cases in which discrete probability distributions
are required. Sometimes it is impossible or very difficult to measure the
life length of a machine on a continuous scale. For example, on-off
switching machines, bulb of photocopier device, etc.

In the last years, many discrete distributions have been derived by
discretizing a known continuous distributions. It are obtained by using the
same method that used to obtain the discrete geometric (DG) distribution
from the continuous exponential distribution. Nakagawa and Osaki (1975)
obtained the discrete Weibull (DW) distribution. A new discrete Weibull
distribution is proposed by Stein and Dattero (1984). Roy (2003) proposed
the discrete normal (DN) distribution. The discrete Rayleigh distribution
(DR) is introduced by Roy (2004). Krishna and Pundir (2009) proposed the
discrete burr (DB) and discrete Pareto (DP) distributions. Gomez and
Calderin (2011) obtained the discrete Lindley\ (DL) distribution. The
discrete generalized exponential (DGE) distribution is proposed by Nekoukhou
et al. (2013). Nekoukhou and Bidram (2015) introduced a new three parameters
distribution and called it the exponentiated discrete Weibull (EDW)
distribution. The CDF and the PMF of the EDW distribution are given
respectively by%
\begin{equation}
F_{EDW}(x;\alpha ,p,\beta )=[1-p^{([x]+1)^{\alpha }}]^{\beta },x\geq 0,
\label{1.1}
\end{equation}%
and 
\begin{eqnarray}
P(X &=&x)=f_{EDW}\left( x;\alpha ,p,\beta \right) =[1-p^{(x+1)^{\alpha
}}]^{\beta }-[1-p^{x^{\alpha }}]^{\beta }  \notag \\
&=&\overset{\infty }{\underset{k=1}{\sum }}(-1)^{k+1}\tbinom{\beta }{k}%
[p^{kx^{\alpha }}-p^{k(x+1)^{\alpha }}],x\in 
\mathbb{N}
_{\circ }=\{0,1,2,...\},  \label{1.3}
\end{eqnarray}%
where $\alpha ,\beta >0$, $0<p<1$ and $[x]$ is the largest integer less than
or equal $x$. For integer $\beta >0,$ the sum in equation (\ref{1.3}) stop
at $\beta .$ There exist some special discrete distributions can be obtained
from EDW distribution as follows:

\begin{enumerate}
\item If $\beta =1,$ then the DW distribution of Nakagawa and Osaki (1975)
is achieved.

\item If $\alpha =1,$ we get the DGE distribution of Nekoukhou et al. (2013).

\item If $\beta =1$ and $\alpha =1,$ then the DG distribution (discrete
exponential (DE) distribution) is obtained.

\item If $\beta =1$ and $\alpha =2,$ then the DR distribution of Roy (2004)
is achieved.

\item If $\alpha =2,$ we get the discrete generalized Rayleigh (DGR)
distribution of Alamatsaz et al. (2016).
\end{enumerate}

It is very useful in simulation study for EDW distribution to know the
following relation: If$\ $the continuous random variable $Y$ has
exponentiated Weibull (EW) distribution, say $Y$\ $\sim EW(\alpha ,\lambda
,\beta ),$ $\lambda =-\ln (p),$ then $X=[Y]\sim EDW(\alpha ,p,\beta ).$ So,
to generate a random sample from the EDW distribution, we first generate a
random sample from a continuous EW distribution by using the inverse CDF
method, and then by considering $X=[Y],$ we find the desired random sample.

On the other hand, the bivariate distributions have been derived and
discussed by many authors which have many applications in the areas such as
engineering, reliability, sports, weather, drought, among others. Until now,
many continuous bivariate distributions based on Marshall and Olkin (1976)
model have been introduced in the literature, see Jose et al. (2009), Kundu
and Gupta (2009), Sarhan et al. (2011), El-Sherpieny et al. (2013), Wagner
and Artur (2013), El- Bassiouny et al. (2016), Rasool and Akbar (2016),
El-Gohary et al. (2016), Mohamed et al. (2017), among others.

Also, many discrete bivariate distributions have been introduced, see
Kocherlakota and Kocherlakota (1992), Kumar (2008), Kemp (2013), Lee and Cha
(2015), Nekoukhou and Kundu (2017), among others.

So, our reasons for introducing the BDEW distribution are the following: to
define a bivariate discrete model having different shapes of the hazard rate
function, and to define a bivariate discrete model having the flexibility
for fitting the real data sets for various phenomena. 

The paper is organized as follows: In Section 2, the BDEW distribution is
defined. Moreover, the joint CDF and the joint PMF are also presented.
Further, some mathematical properties of the BDEW distribution such as the
joint PGF, the marginal CDF, the marginal PMF, the conditional PMF of $X_{1}$
given $X_{2}=x_{2}$, the conditional CDF of $X_{1}$ given $X_{2}\leq x_{2},$
the conditional CDF of $X_{1}$ given $X_{2}=x_{2}$, the conditional
expectation of $X_{1}$ given $X_{2}=x_{2}$ and some other results are
presented in Section 3. In Section 4, some reliability studies are
introduced. In Section 5, the parameters of the BDEW distribution are
estimated by the maximum likelihood method. In Section 6, two real data sets
are analyzed to show the importance of the proposed distribution. Finally,
Section 7 offers some concluding remarks.

\section{\textbf{The BDEW Distribution }}

Suppose that $V_{i},$ $i=1,2,3$ are three independently distributed random
variables, and let $V_{i}\sim EDW(\alpha ,p,\beta _{i})$. If $X_{1}=\max
\{V_{1},V_{3}\}$ and $X_{2}=\max \{V_{2},V_{3}\}$, \ then the bivariate
vector $(X_{1},X_{2})$ has the BDEW distribution with the parameter vector $%
\mathbf{\Omega }=(\alpha ,p,\beta _{1},\beta _{2},\beta _{3})^{T}$.

\paragraph{\textit{Lemma 1:}}

If $(X_{1},X_{2})\thicksim $BDEW($\mathbf{\Omega }$), then the joint CDF of $%
(X_{1},X_{2})$ is given by

\begin{eqnarray}
F_{X_{1},X_{2}}(x_{1},x_{2}) &=&[1-p^{(x_{1}+1)^{\alpha }}]^{\beta
_{1}}[1-p^{(x_{2}+1)^{\alpha }}]^{\beta _{2}}[1-p^{(z+1)^{\alpha }}]^{\beta
_{3}}  \notag \\[0.02in]
&=&\left\{ 
\begin{array}{c}
F_{1}(x_{1},x_{2})\text{ \ \ \ \ \ \ \ \ \ \ \ \ \ \ \ \ \ \ \ \ \ \ \ \ \ \
\ \ \ \ if \ \ }x_{1}<x_{2} \\ 
F_{2}(x_{1},x_{2})\text{ \ \ \ \ \ \ \ \ \ \ \ \ \ \ \ \ \ \ \ \ \ \ \ \ \ \
\ \ \ \ \ if \ \ }x_{2}<x_{1} \\ 
F_{3}(x)\text{\ \ \ \ \ \ \ \ \ \ \ \ \ \ \ \ \ \ \ \ \ \ \ \ \ \ \ \ \ \ \
if \ \ }x_{1}=x_{2}=x,%
\end{array}%
\right.  \label{5}
\end{eqnarray}%
where $x_{1},x_{2}\in 
\mathbb{N}
_{\circ }$, $z=\min \{x_{1},x_{2}\}$ and $F_{1}(x_{1},x_{2}),$ $%
F_{2}(x_{1},x_{2}),$ $F_{3}(x)$ are given by

\begin{equation*}
F_{1}(x_{1},x_{2})=F_{EDW}(x_{1};\alpha ,p,\beta _{1}+\beta
_{3})F_{EDW}(x_{2};\alpha ,p,\beta _{2}),
\end{equation*}

\begin{equation*}
F_{2}(x_{1},x_{2})=F_{EDW}(x_{1};\alpha ,p,\beta _{1})F_{EDW}(x_{2};\alpha
,p,\beta _{2}+\beta _{3})\text{ }
\end{equation*}%
and

\begin{equation*}
F_{3}(x)=F_{EDW}(x;\alpha ,p,\beta _{1}+\beta _{2}+\beta _{3}).
\end{equation*}

\paragraph{\textit{Proof:}}

The joint CDF of the random variables $X_{1}$ and $X_{2}$ is defined as
follows 

\begin{eqnarray*}
F_{X_{1},X_{2}}(x_{1},x_{2}) &=&P(X_{1}\leq x_{1},X_{2}\leq x_{2}) \\
&=&P(\max \{V_{1},V_{3}\}\leq x_{1},\max \{V_{2},V_{3}\}\leq x_{2}) \\
&=&P(V_{1}\leq x_{1},V_{2}\leq x_{2},V_{3}\leq \min \{x_{1},x_{2}\}).
\end{eqnarray*}%
Since, the random variables $V_{i},i=1,2,3$ are independent, we obtain

\begin{eqnarray}
F_{X_{1},X_{2}}(x_{1},x_{2}) &=&P(V_{1}\leq x_{1})P(V_{2}\leq
x_{2})P(V_{3}\leq \min \{x_{1},x_{2}\})  \notag \\
&=&F_{EDW}(x_{1};\alpha ,p,\beta _{1})F_{EDW}(x_{2};\alpha ,p,\beta
_{2})F_{EDW}(z;\alpha ,p,\beta _{3}).  \label{1.9}
\end{eqnarray}%
By substituting from (\ref{1.1}) into (\ref{1.9}), we get (\ref{5}), which
complete the proof.

\subsection{The joint PMF}

The joint PMF of the bivariate vector $(X_{1},X_{2})$ can be easily obtained
by using the following relation:%
\begin{equation}
f_{X_{1},X_{2}}(x_{1},x_{2})=F_{X_{1},X_{2}}(x_{1},x_{2})-F_{X_{1},X_{2}}(x_{1}-1,x_{2})-F_{X_{1},X_{2}}(x_{1},x_{2}-1)+F_{X_{1},X_{2}}(x_{1}-1,x_{2}-1).
\end{equation}%
The joint PMF of $(X_{1},X_{2})$ for $x_{1},x_{2}\in 
\mathbb{N}
_{\circ }$ is given by

\begin{equation}
f_{X_{1},X_{2}}(x_{1},x_{2})=\left\{ 
\begin{array}{c}
f_{1}(x_{1},x_{2})\text{ \ \ \ \ \ \ \ \ \ \ \ \ \ if \ \ }x_{1}<x_{2} \\ 
f_{2}(x_{1},x_{2})\text{ \ \ \ \ \ \ \ \ \ \ \ \ \ if \ \ }x_{2}<x_{1} \\ 
f_{3}(x)\text{ \ \ \ \ \ \ \ \ \ \ \ \ \ if \ \ }x_{1}=x_{2}=x,%
\end{array}%
\right.  \label{1.10}
\end{equation}%
where

\begin{eqnarray*}
f_{1}(x_{1},x_{2})\text{ \ } &=&\left( [1-p^{(x_{1}+1)^{\alpha }}]^{\beta
_{1}+\beta _{3}}-[1-p^{x_{1}^{\alpha }}]^{\beta _{1}+\beta _{3}}\right)
\left( [1-p^{(x_{2}+1)^{\alpha }}]^{\beta _{2}}-[1-p^{x_{2}^{\alpha
}}]^{\beta _{2}}\right) \\
&=&f_{EDW}(x_{1};\alpha ,p,\beta _{1}+\beta _{3})f_{EDW}(x_{2};\alpha
,p,\beta _{2}),
\end{eqnarray*}

\begin{eqnarray*}
f_{2}(x_{1},x_{2})\text{ } &=&\left( [1-p^{(x_{1}+1)^{\alpha }}]^{\beta
_{1}}-[1-p^{x_{1}^{\alpha }}]^{\beta _{1}}\right) \left(
[1-p^{(x_{2}+1)^{\alpha }}]^{\beta _{2}+\beta _{3}}-[1-p^{x_{2}^{\alpha
}}]^{\beta _{2}+\beta _{3}}\right) \\
&=&f_{EDW}(x_{1};\alpha ,p,\beta _{1})f_{EDW}(x_{2};\alpha ,p,\beta
_{2}+\beta _{3})
\end{eqnarray*}%
and

\begin{eqnarray*}
f_{3}(x)\text{ } &=&p_{1}\left( [1-p^{(x+1)^{\alpha }}]^{\beta _{2}+\beta
_{3}}-[1-p^{x^{\alpha }}]^{\beta _{2}+\beta _{3}}\right) -p_{2}\left(
[1-p^{(x+1)^{\alpha }}]^{\beta _{2}}-[1-p^{x^{\alpha }}]^{\beta _{2}}\right)
\\
&=&p_{1}f_{EDW}(x;\alpha ,p,\beta _{2}+\beta _{3})-p_{2}f_{EDW}(x;\alpha
,p,\beta _{2}),
\end{eqnarray*}%
where%
\begin{equation*}
p_{1}=[1-p^{(x+1)^{\alpha }}]^{\beta _{1}}\text{, }p_{2}=[1-p^{x^{\alpha
}}]^{\beta _{1}+\beta _{3}}.
\end{equation*}

The scatter plot of the joint PMF for the BDEW distribution is given in
Figure 1. As expected, the joint PMF for the BDEW distribution can take
varies shapes depending on the values of its parameter vector $\mathbf{%
\Omega }.$ And so, this distribution is more\ flexible to provide a better
fit to variety of data sets.

\begin{eqnarray*}
&&\FRAME{itbpFU}{2.8383in}{2.8383in}{0in}{\Qcb{{}}}{}{pmf1.eps}{\special%
{language "Scientific Word";type "GRAPHIC";maintain-aspect-ratio
TRUE;display "USEDEF";valid_file "F";width 2.8383in;height 2.8383in;depth
0in;original-width 5.7519in;original-height 5.7519in;cropleft "0";croptop
"1";cropright "1";cropbottom "0";filename '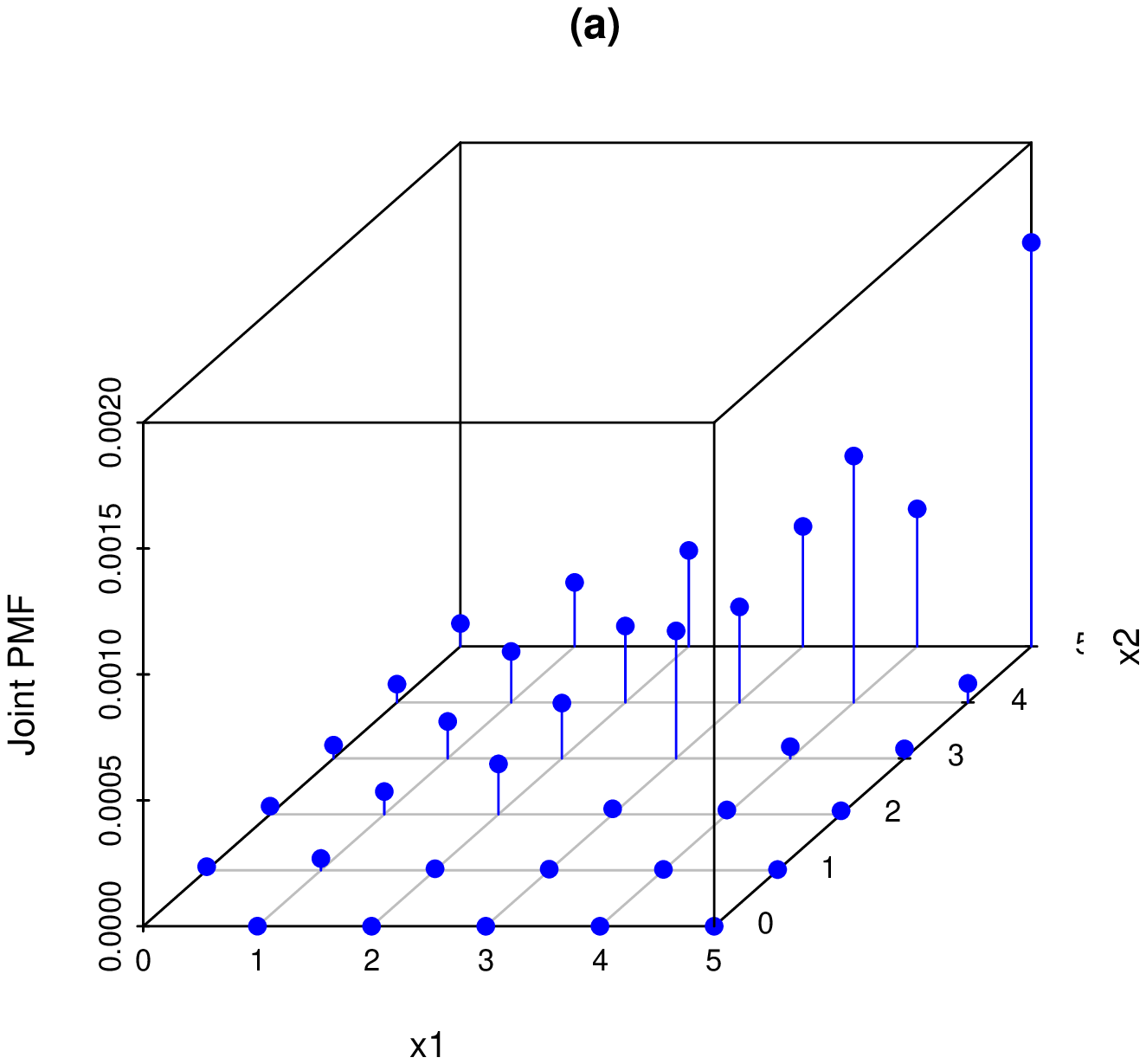';file-properties
"XNPEU";}}\FRAME{itbpFU}{2.8383in}{2.8383in}{0in}{\Qcb{{}}}{}{pmf2.eps}{%
\special{language "Scientific Word";type "GRAPHIC";maintain-aspect-ratio
TRUE;display "USEDEF";valid_file "F";width 2.8383in;height 2.8383in;depth
0in;original-width 5.7519in;original-height 5.7519in;cropleft "0";croptop
"1";cropright "1";cropbottom "0";filename '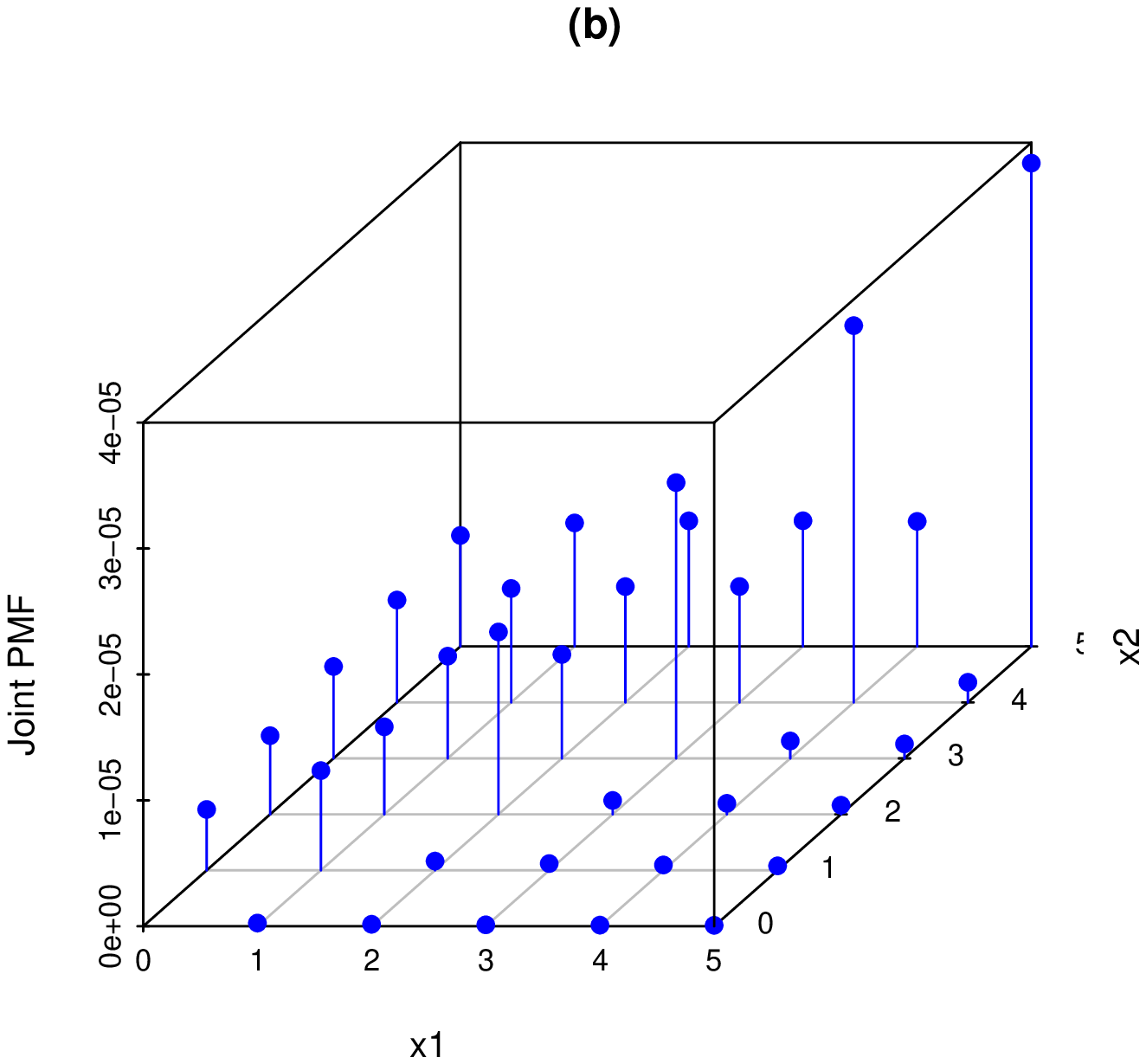';file-properties
"XNPEU";}} \\
&&\FRAME{itbpFU}{2.8383in}{2.8383in}{0in}{\Qcb{{}}}{}{pmf3.eps}{\special%
{language "Scientific Word";type "GRAPHIC";maintain-aspect-ratio
TRUE;display "USEDEF";valid_file "F";width 2.8383in;height 2.8383in;depth
0in;original-width 5.7519in;original-height 5.7519in;cropleft "0";croptop
"1";cropright "1";cropbottom "0";filename '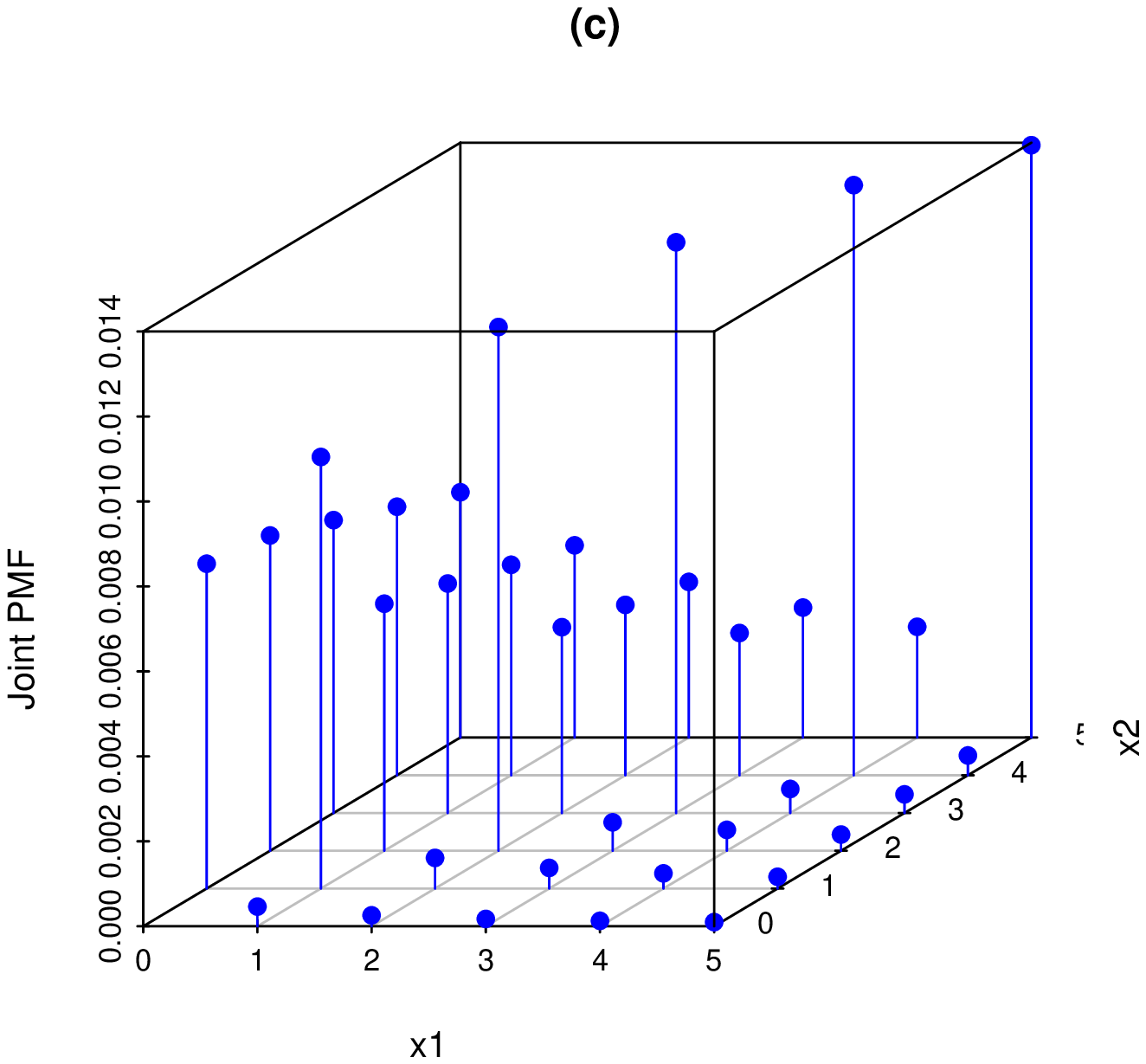';file-properties
"XNPEU";}}\FRAME{itbpFU}{2.8383in}{2.8383in}{0in}{\Qcb{{}}}{}{pmf4.eps}{%
\special{language "Scientific Word";type "GRAPHIC";maintain-aspect-ratio
TRUE;display "USEDEF";valid_file "F";width 2.8383in;height 2.8383in;depth
0in;original-width 5.7519in;original-height 5.7519in;cropleft "0";croptop
"1";cropright "1";cropbottom "0";filename '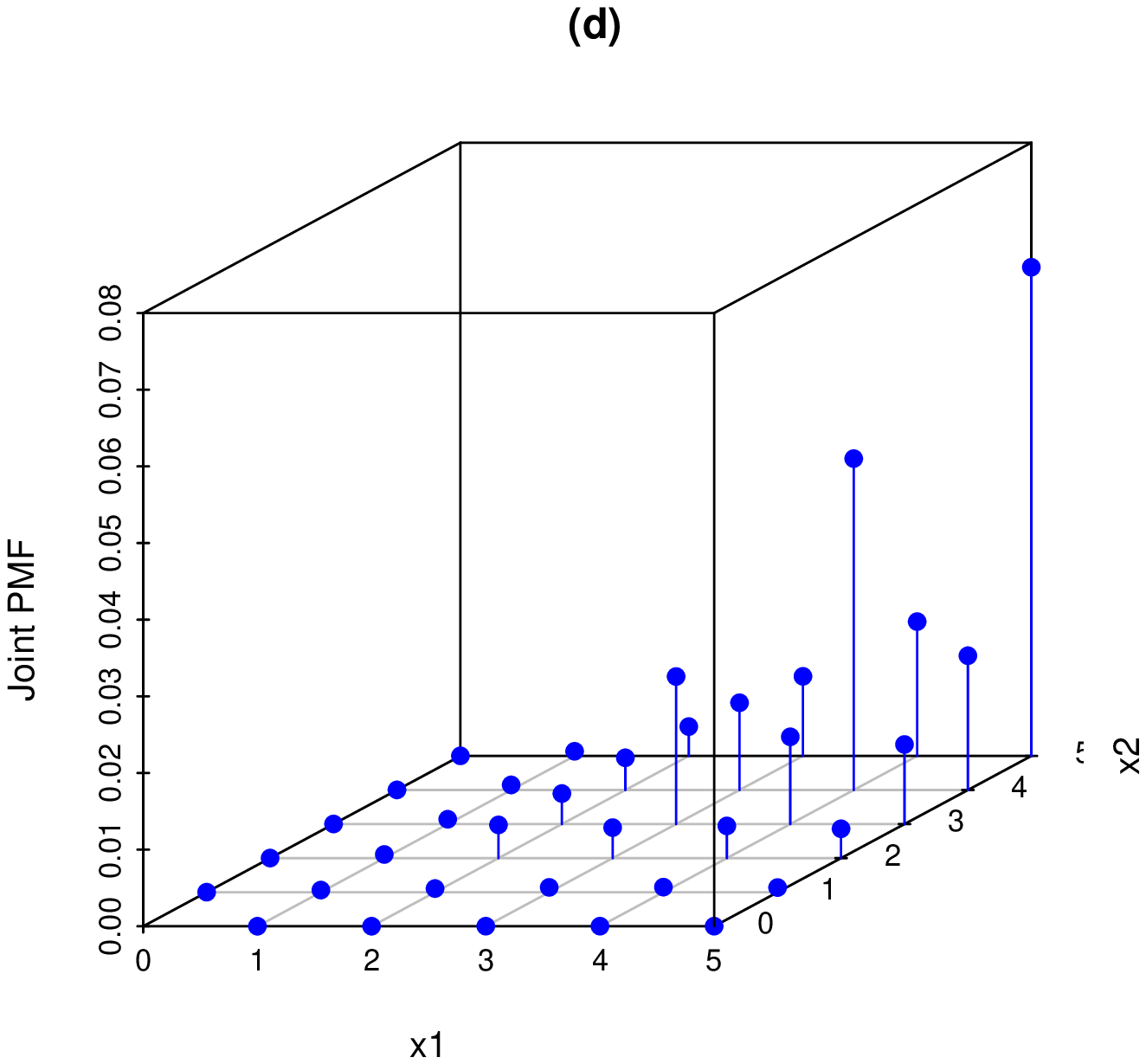';file-properties
"XNPEU";}}
\end{eqnarray*}

\begin{center}
\textbf{Figure 1.} Scatter plot of the joint PMF of the BDEW distribution
for different values of its parameter vector $\mathbf{\Omega }=(\alpha
,p,\beta _{1},\beta _{2},\beta _{3})$ : \textbf{(a)} $\mathbf{\Omega }%
=(0.9,0.9,0.3,3.3,2),$ \textbf{(b)} $\mathbf{\Omega }=(0.5,0.9,0.3,3.3,2),$ 
\textbf{(c)} $\mathbf{\Omega }=(0.5,0.6,0.3,3.3,2)$ and \textbf{(d)} $%
\mathbf{\Omega }=(1.5,0.9,2,2,2)$.
\end{center}

\subsubsection{Special case}

Some special bivariate discrete distributions are achieved from the BDEW
distribution as follows:

\begin{enumerate}
\item If $\alpha =1$, the bivariate discrete generalized exponential (BDGE)
distribution\ of Nekoukhou and kundu (2017) is obtained.

\item If $\alpha =2$, we get the bivariate discrete generalized Rayleigh
(BDGR) distribution.

\item If $\alpha =1,$ $\beta _{3}=\beta ,$ $0<\beta <1$ and also $\beta
_{1}=\beta _{2}=1-\beta ,$ then we have a new bivariate geomatric (NBG)
distribution with two parameters $0<\alpha <1$ and $0<p<1,$ see, Nekoukhou
and kundu (2017). The joint CDF of the NBG distribution is%
\begin{equation*}
F_{X_{1},X_{2}}(x_{1},x_{2})=[1-p^{x_{1}+1}]^{1-\beta
}[1-p^{x_{2}+1}]^{1-\beta }[1-p^{z+1}]^{\beta },z=\min \{x_{1},x_{2}\}.
\end{equation*}
\end{enumerate}

\section{Statistical Properties}

\subsection{The joint probability generating function}

The joint probability generating function (PGF) for any bivariate
distribution is very useful and important, because we can use it to find the
varies moments, and also the product moments as infinite series. The PGF of
the BDEW distribution is mentioned in the following theorem.

\paragraph{\textit{Theorem 1:}}

If $(X_{1},X_{2})\thicksim $BDEW($\mathbf{\Omega }$), then the PGF of $%
(X_{1},X_{2})$ is given by

\begin{eqnarray}
G(u,v) &=&\underset{i=0}{\overset{\infty }{\sum }}\overset{\infty }{\underset%
{j=i+1}{\sum }}\underset{k=1}{\overset{\infty }{\sum }}\overset{\infty }{%
\underset{l=1}{\sum }}(-1)^{k+l}\binom{\beta _{1}+\beta _{3}}{k}\binom{\beta
_{2}}{l}[p^{ki^{\alpha }}-p^{k(i+1)^{\alpha }}][p^{kj^{\alpha
}}-p^{k(j+1)^{\alpha }}]u^{i}v^{j}  \notag \\
&&+\underset{j=0}{\overset{\infty }{\sum }}\overset{\infty }{\underset{i=j+1}%
{\sum }}\underset{k=1}{\overset{\infty }{\sum }}\overset{\infty }{\underset{%
l=1}{\sum }}(-1)^{k+l}\binom{\beta _{1}}{k}\binom{\beta _{2}+\beta _{3}}{l}%
[p^{ki^{\alpha }}-p^{k(i+1)^{\alpha }}][p^{lj^{\alpha }}-p^{l(j+1)^{\alpha
}}]u^{i}v^{j}  \notag \\
&&+\underset{j=0}{\overset{\infty }{\sum }}\underset{i=0}{\overset{\infty }{%
\sum }}\underset{k=1}{\overset{\infty }{\sum }}(-1)^{j+k+l}\binom{\beta _{1}%
}{k}\binom{\beta _{2}+\beta _{3}}{j}p^{j(i+1)^{\alpha }}[p^{ki^{\alpha
}}-p^{k(i+1)^{\alpha }}]u^{i}v^{i}  \notag \\
&&-\underset{j=0}{\overset{\infty }{\sum }}\underset{i=0}{\overset{\infty }{%
\sum }}\underset{k=1}{\overset{\infty }{\sum }}(-1)^{j+k+l}\binom{\beta
_{1}+\beta _{3}}{k}\binom{\beta _{2}}{j}p^{ji^{\alpha }}[p^{ki^{\alpha
}}-p^{k(i+1)^{\alpha }}]u^{i}v^{i},  \label{1.14}
\end{eqnarray}%
where $\left\vert V\right\vert <1$ and $\left\vert v\right\vert <1.$

\paragraph{\textit{Proof:}}

From the definition of the joint PGF of $(X_{1},X_{2}),$ we get

\begin{equation}
G(u,v)=E(u^{X_{1}}v^{X_{2}})=\overset{\infty }{\underset{i=0}{\sum }}%
\underset{j=0}{\overset{\infty }{\sum }}%
P(X_{1}=i,X_{2}=j)u^{i}v^{j}=I+II+III,  \label{1.15}
\end{equation}%
where

\begin{equation}
I=\underset{i=0}{\overset{\infty }{\sum }}\overset{\infty }{\underset{j=i+1}{%
\sum }}P(X_{1}=i,X_{2}=j)u^{i}v^{j},  \label{1.16}
\end{equation}

\begin{equation}
II=\underset{j=0}{\overset{\infty }{\sum }}\overset{\infty }{\underset{i=j+1}%
{\sum }}P(X_{1}=i,X_{2}=j)u^{i}v^{j}  \label{1.17}
\end{equation}%
and

\begin{equation}
III=\overset{\infty }{\underset{i=0}{\sum }}\underset{j=0}{\overset{\infty }{%
\sum }}P(X_{1}=i,X_{2}=i)u^{i}v^{i}.  \label{1.18}
\end{equation}%
By substituting from (\ref{1.10}) into (\ref{1.16}), (\ref{1.17}) and (\ref%
{1.18}), we find

\begin{eqnarray}
I &=&\underset{i=0}{\overset{\infty }{\sum }}\overset{\infty }{\underset{%
j=i+1}{\sum }}f_{EDW}(i;\alpha ,p,\beta _{1}+\beta _{3})f_{EDW}(j;\alpha
,p,\beta _{2})u^{i}v^{j}  \notag \\
&=&\underset{i=0}{\overset{\infty }{\sum }}\overset{\infty }{\underset{j=i+1}%
{\sum }}\underset{k=1}{\overset{\infty }{\sum }}\overset{\infty }{\underset{%
l=1}{\sum }}(-1)^{k+l}\binom{\beta _{1}+\beta _{3}}{k}\binom{\beta _{2}}{l}%
\left[ p^{ki^{\alpha }}-p^{k(i+1)^{\alpha }}\right] \left[ p^{lj^{\alpha
}}-p^{l(j+1)^{\alpha }}\right] u^{i}v^{j},  \notag \\
&&  \label{1.19}
\end{eqnarray}

\begin{eqnarray}
II &=&\underset{j=0}{\overset{\infty }{\sum }}\overset{\infty }{\underset{%
i=j+1}{\sum }}f_{EDW}(i;\alpha ,p,\beta _{1})f_{EDW}(j;\alpha ,p,\beta
_{2}+\beta _{3})u^{i}v^{j}  \notag \\
&=&\underset{j=0}{\overset{\infty }{\sum }}\overset{\infty }{\underset{i=j+1}%
{\sum }}\underset{k=1}{\overset{\infty }{\sum }}\overset{\infty }{\underset{%
l=1}{\sum }}(-1)^{k+l}\binom{\beta _{1}}{k}\binom{\beta _{2}+\beta _{3}}{l}%
\left[ p^{ki^{\alpha }}-p^{k(i+1)^{\alpha }}\right] \left[ p^{lj^{\alpha
}}-p^{l(j+1)^{\alpha }}\right] u^{i}v^{j}  \notag \\
&&  \label{1.20}
\end{eqnarray}%
and

\begin{eqnarray}
III &=&\overset{\infty }{\underset{i=0}{\sum }}\underset{j=0}{\overset{%
\infty }{\sum }}\left[ 1-p^{(i+1)^{\alpha }}\right] ^{\beta
_{1}}f_{EDW}(i;\alpha ,p,\beta _{2}+\beta _{3})u^{i}v^{i}  \notag \\
&&-\overset{\infty }{\underset{i=0}{\sum }}\underset{j=0}{\overset{\infty }{%
\sum }}\left[ 1-p^{i^{\alpha }}\right] ^{\beta _{1}+\beta
_{3}}f_{EDW}(i;\alpha ,p,\beta _{2})u^{i}v^{i}  \notag \\
&=&\underset{j=0}{\overset{\infty }{\sum }}\underset{i=0}{\overset{\infty }{%
\sum }}\underset{k=1}{\overset{\infty }{\sum }}(-1)^{j+k+l}\binom{\beta _{1}%
}{k}\binom{\beta _{2}+\beta _{3}}{j}p^{j(i+1)^{\alpha }}[p^{ki^{\alpha
}}-p^{k(i+1)^{\alpha }}]u^{i}v^{i}  \notag \\
&&-\underset{j=0}{\overset{\infty }{\sum }}\underset{i=0}{\overset{\infty }{%
\sum }}\underset{k=1}{\overset{\infty }{\sum }}(-1)^{j+k+l}\binom{\beta
_{1}+\beta _{3}}{k}\binom{\beta _{2}}{j}p^{ji^{\alpha }}\left[ p^{ki^{\alpha
}}-p^{k(i+1)^{\alpha }}\right] u^{i}v^{i}.  \label{1.21}
\end{eqnarray}%
Substituting from (\ref{1.19}), (\ref{1.20}) and (\ref{1.21}) into (\ref%
{1.15}), we get (\ref{1.14}), which complete the proof.

\subsection{The marginal CDF and PMF of $X_{1}$ and $X_{2}$}

\paragraph{\textit{Lemma 2:}}

The marginal CDF of $X_{i}$ $,(i=1,2)$ is given by

\begin{equation}
F_{X_{i}}(x_{i})=F_{EDW}(x_{i};\alpha ,p,\beta _{i}+\beta _{3})=\left[
1-p^{(x_{i}+1)^{\alpha }}\right] ^{\beta _{i}+\beta _{3}},x_{i}\in 
\mathbb{N}
_{\circ }.  \label{1.222}
\end{equation}

\paragraph{\textit{Proof:}}

The CDF of $X_{i}$ given by%
\begin{equation*}
F_{X_{i}}(x_{i})=P(X_{i}\leq x_{i})=P(\max \{V_{i},V_{3}\}\leq
x_{i})=P(V_{i}\leq x_{i},V_{3}\leq x_{i}).
\end{equation*}%
Because the random variables $V_{i},(i=1,2)$ and $V_{3}$ are independent, we
obtain

\begin{eqnarray*}
F_{X_{i}}(x_{i}) &=&P(V_{i}\leq x_{i})P(V_{3}\leq x_{i}) \\
&=&\left[ 1-p^{(x_{i}+1)^{\alpha }}\right] ^{\beta _{i}}\left[
1-p^{(x_{i}+1)^{\alpha }}\right] ^{\beta _{3}} \\
&=&\left[ 1-p^{(x_{i}+1)^{\alpha }}\right] ^{\beta _{i}+\beta
_{3}}=F_{EDW}(x_{i};\alpha ,p,\beta _{i}+\beta _{3}).
\end{eqnarray*}

\paragraph{\textit{Remark:}}

The marginal PMF of $X_{i},(i=1,2),$ corresponding to (\ref{1.222}) is%
\begin{eqnarray}
f_{X_{i}}(x_{i}) &=&f_{EDW}(x_{i};\alpha ,p,\beta _{i}+\beta _{3})  \notag \\
&=&[1-p^{(x_{i}+1)^{\alpha }}]^{\beta _{i}+\beta _{3}}-[1-p^{x_{i}^{\alpha
}}]^{\beta _{i}+\beta _{3}},x_{i}\in 
\mathbb{N}
_{\circ }.  \label{1.233}
\end{eqnarray}

\subsection{The conditional PMF of $X_{1}\ $given $X_{2}=x_{2}$}

\bigskip The conditional PMF of $(X_{1}\mid X_{2}=x_{2}),$ say $f_{X_{1}\mid
X_{2}=x_{2}}(x_{1}\mid x_{2}),$ is given by%
\begin{equation}
f_{X_{1}\mid X_{2}=x_{2}}(x_{1}\mid x_{2})=\left\{ 
\begin{array}{c}
f_{1}(x_{1}\mid x_{2})\text{ \ \ \ \ \ \ \ \ \ \ \ \ \ \ \ \ \ \ \ if \ \ }%
0\leq x_{1}<x_{2} \\ 
f_{2}(x_{1}\mid x_{2})\text{ \ \ \ \ \ \ \ \ \ \ \ \ \ \ \ \ \ \ \ if \ \ }%
0\leq x_{2}<x_{1} \\ 
f_{3}(x_{1}\mid x_{2})\text{ \ \ \ \ \ \ \ \ \ \ if \ \ }0\leq x_{1}=x_{2}=x,%
\text{\ }%
\end{array}%
\right.  \label{1.2222}
\end{equation}%
where%
\begin{eqnarray*}
f_{1}(x_{1} &\mid &x_{2})=\frac{\left( [1-p^{(x_{1}+1)^{\alpha }}]^{\beta
_{1}+\beta _{3}}-[1-p^{x_{1}^{\alpha }}]^{\beta _{1}+\beta _{3}}\right)
\left( [1-p^{(x_{2}+1)^{\alpha }}]^{\beta _{2}}-[1-p^{x_{2}^{\alpha
}}]^{\beta _{2}}\right) }{[1-p^{(x_{2}+1)^{\alpha }}]^{\beta _{2}+\beta
_{3}}-[1-p^{x_{2}^{\alpha }}]^{\beta _{2}+\beta _{3}}}, \\
f_{2}(x_{1} &\mid &x_{2})=[1-p^{(x_{1}+1)^{\alpha }}]^{\beta
_{1}}-[1-p^{x_{1}^{\alpha }}]^{\beta _{1}}
\end{eqnarray*}%
and%
\begin{equation*}
f_{3}(x_{1}\mid x_{2})=[1-p^{(x+1)^{\alpha }}]^{\beta _{1}}-\frac{%
[1-p^{x^{\alpha }}]^{\beta _{1}+\beta _{3}}\left( [1-p^{(x+1)^{\alpha
}}]^{\beta _{2}}-[1-p^{x^{\alpha }}]^{\beta _{2}}\right) }{%
[1-p^{(x+1)^{\alpha }}]^{\beta _{2}+\beta _{3}}-[1-p^{x^{\alpha }}]^{\beta
_{2}+\beta _{3}}}.
\end{equation*}%
Equation (\ref{1.2222}) can be getting by using the following relation:

\begin{equation*}
f_{X_{1}\mid X_{2}=x_{2}}(x_{1}\mid x_{2})=\frac{P(X_{1}=x_{1},X_{2}=x_{2})}{%
P(X_{2}=x_{2})}.
\end{equation*}

\subsection{The conditional CDF of $X_{1}$ given $X_{2}\leq x_{2}$}

The conditional CDF of $(X_{1}\mid X_{2}\leq x_{2}),$ say $F_{X_{1}\mid
X_{2}\leq x_{2}}(x_{1}),$ is given by

\begin{equation}
F_{X_{1}\mid X_{2}\leq x_{2}}(x_{1})=\left\{ 
\begin{array}{c}
\lbrack 1-p^{(x_{1}+1)^{\alpha }}]^{\beta _{1}+\beta
_{3}}[1-p^{(x_{2}+1)^{-\beta _{3}}}]\text{\ \ \ \ \ if \ \ }0\leq x_{1}<x_{2}
\\ 
\lbrack 1-p^{(x_{1}+1)^{\alpha }}]^{\beta _{1}}\text{ \ \ \ \ \ \ \ \ \ \ \
\ \ \ \ \ \ \ \ \ \ \ \ \ \ \ \ \ \ \ if \ \ }0\leq x_{2}<x_{1} \\ 
\lbrack 1-p^{(x+1)^{\alpha }}]^{\beta _{1}}\text{ \ \ \ \ \ \ \ \ \ \ \ \ \
\ \ \ \ \ \ \ \ \ \ if \ \ }0\leq x_{1}=x_{2}=x\text{.}%
\end{array}%
\right. \text{ \ \ \ \ }  \label{1.22}
\end{equation}%
Equation (\ref{1.22}) can be getting by using the following relation%
\begin{equation*}
F_{X_{1}\mid X_{2}\leq x_{2}}(x_{1})=\frac{P(X_{1}\leq x_{1},X_{2}\leq x_{2})%
}{P(X_{2}\leq x_{2})}.
\end{equation*}

\subsection{The conditional CDF of $X_{1}$ given $X_{2}=x_{2}$}

The conditional CDF of $(X_{1}\mid X_{2}=x_{2}),$ say $F_{X_{1}\mid
X_{2}=x_{2}}(x_{1}),$ is given by%
\begin{equation}
F_{X_{1}\mid X_{2}=x_{2}}(x_{1})=\left\{ 
\begin{array}{c}
F_{1}(x_{1}\mid x_{2})\text{\ \ \ \ \ \ \ \ \ \ \ \ if \ \ }0\leq x_{1}<x_{2}%
\text{\ } \\ 
F_{2}(x_{1}\mid x_{2})\text{\ \ \ \ \ \ \ \ \ \ \ \ \ \ if \ \ }0\leq
x_{2}<x_{1} \\ 
F_{3}(x_{1}\mid x_{2})\text{\ \ \ \ \ if \ \ }0\leq x_{1}=x_{2}=x,\text{\ }%
\end{array}%
\right.  \label{1.23}
\end{equation}%
where

\begin{equation*}
F_{1}(x_{1}\mid x_{2})=\frac{F_{EDW}(x_{1};\alpha ,p,\beta _{1}+\beta _{3}%
\text{ })f_{EDW}(x_{2};\alpha ,p,\beta _{2})}{f_{EDW}(x_{2};\alpha ,p,\beta
_{2}+\beta _{3})},
\end{equation*}

\begin{equation*}
F_{2}(x_{1}\mid x_{2})=F_{EDW}(x_{1};\alpha ,p,\beta _{1})
\end{equation*}%
and

\begin{equation*}
F_{3}(x_{1}\mid x_{2})=\frac{F_{EDW}(x;\alpha ,p,\beta _{1}+\beta _{2}+\beta
_{3})-\text{ }F_{EDW}(x;\alpha ,p,\beta _{1})\text{ }F_{EDW}(x-1;\alpha
,p,\beta _{2}+\beta _{3}\text{ })\text{\ }}{f_{EDW}(x_{2};\alpha ,p,\beta
_{2}+\beta _{3})}.
\end{equation*}%
Equation (\ref{1.23}) can be getting by using the following relation

\begin{equation*}
F_{X_{1}\mid X_{2}=x_{2}}(x_{1})=\frac{P(X_{1}\leq x_{1},X_{2}=x_{2})}{%
P(X_{2}=x_{2})},
\end{equation*}%
which

\begin{equation*}
F_{1}(x_{1}\mid x_{2})=\frac{\overset{x_{1}}{\underset{j=0}{\dsum }}%
P(X_{1}=j,X_{2}=x_{2})}{P(X_{2}=x_{2})},
\end{equation*}

\begin{equation*}
F_{2}(x_{1}\mid x_{2})=\frac{\overset{x_{1}-1}{\underset{j=0}{\dsum }}%
P(X_{1}=j,X_{2}=x_{2})+P(X_{1}=x_{2},X_{2}=x_{2})+\overset{x_{2}-1}{\underset%
{j=x_{1}+1}{\dsum }}P(X_{1}=j,X_{2}=x_{2})}{P(X_{2}=x_{2})}
\end{equation*}%
and

\begin{equation*}
F_{3}(x_{1}\mid x_{2})=\frac{\overset{x-1}{\underset{j=0}{\dsum }}%
P(X_{1}=j,X_{2}=x_{2})+P(X_{1}=x,X_{2}=x)}{P(X_{2}=x_{2})}.
\end{equation*}

\subsection{The conditional expectation of $X_{1}\ $given $X_{2}=x_{2}$}

\paragraph{\textit{Lemma 3:}}

The conditional expectation of $(X_{1}\mid X_{2}=x_{2}),$ say $E(X_{1}\mid
X_{2}=x_{2}),$ is given by%
\begin{eqnarray}
E(X_{1} &\mid &X_{2}=x_{2})=\frac{[1-p^{(x_{2}+1)^{\alpha }}]^{\beta
_{2}}-[1-p^{x_{2}^{\alpha }}]^{\beta _{2}}}{[1-p^{(x_{2}+1)^{\alpha
}}]^{\beta _{2}+\beta _{3}}-[1-p^{x_{2}^{\alpha }}]^{\beta _{2}+\beta _{3}}}
\notag \\
&&\times \overset{\infty }{\underset{x_{1}=x_{2}+1}{\sum }}x_{1}\left(
[1-p^{(x_{1}+1)^{\alpha }}]^{\beta _{1}+\beta _{3}}-[1-p^{x_{1}^{\alpha
}}]^{\beta _{1}+\beta _{3}}\right)  \notag \\
&&+\overset{x_{2}-1}{\underset{x_{1}=0}{\sum }}x_{1}\left(
[1-p^{(x_{1}+1)^{\alpha }}]^{\beta _{1}}-[1-p^{x_{1}^{\alpha }}]^{\beta
_{1}}\right) +x_{2}[1-p^{(x_{2}+1)^{\alpha }}]^{\beta _{1}}  \notag \\
&&-\frac{x_{2}[1-p^{x_{2}^{\alpha }}]^{\beta _{1}+\beta _{3}}\left(
[1-p^{(x_{2}+1)^{\alpha }}]^{\beta _{2}}-[1-p^{x_{2}^{\alpha }}]^{\beta
_{2}}\right) }{[1-p^{(x_{2}+1)^{\alpha }}]^{\beta _{2}+\beta
_{3}}-[1-p^{x_{2}^{\alpha }}]^{\beta _{2}+\beta _{3}}}.  \label{1.24}
\end{eqnarray}

\paragraph{\textit{Proof:}}

\begin{eqnarray}
E(X_{1} &\mid &X_{2}=x_{2})=\overset{\infty }{\underset{x_{1}=0}{\sum }}%
x_{1}f_{X_{1}\mid X_{2}=x_{2}}(x_{1}\mid x_{2})  \notag \\
&=&\overset{\infty }{\underset{x_{1}=x_{2}+1}{\sum }}x_{1}f_{1}(x_{1}\mid
x_{2})+\overset{x_{2}-1}{\underset{x_{1}=0}{\sum }}x_{1}f_{2}(x_{1}\mid
x_{2})+x_{2}f_{3}(x_{1}\mid x_{2}).  \label{1.25}
\end{eqnarray}%
Substituting from (\ref{1.2222}) into (\ref{1.25}), the equation (\ref{1.24}%
) is obtained, which complete the proof.

\subsection{Other results}

\paragraph{\textit{Result 1:}}

If $(Z_{1},Z_{2})\ $has a bivariate continuous exponentiated Weibull
distribution with parameters $\alpha ,\lambda ,\beta _{1},\beta _{2}$ and $%
\beta _{3}$, then $(X_{1},X_{2})$ $\sim BDEW(\mathbf{\Omega })$, where $%
X_{i}=\left[ Z_{i}\right] ,i=1,2$, $p=e^{-\lambda }$ and $\left[ Z_{i}\right]
$ is the largest integer less than or equal $Z_{i}$.

\paragraph{\textit{Proof:}}

Obvious.

\paragraph{\textit{Result 2:}}

The stress-strenght probability of $(X_{1},X_{2})$ $\sim BDEW(\alpha
,p,\beta _{1},\beta _{2},\beta _{3})$ is given by

\begin{equation}
P(X_{1}<X_{2})=\overset{\infty }{\underset{i=0}{\sum }}\{[1-p^{(i+2)^{\alpha
}}]^{\beta _{2}}-[1-p^{(i+1)^{\alpha }}]^{\beta _{2}}\}[1-p^{(i+1)^{\alpha
}}]^{\beta _{1}+\beta _{3}}.  \label{1.26}
\end{equation}

\paragraph{\textit{Proof:}}

\begin{eqnarray*}
P(X_{1} &<&X_{2})=\overset{\infty }{\underset{i=0}{\sum }}\overset{j}{%
\underset{i=0}{\sum }}P(X_{2}=j+1,X_{1}=i) \\
&=&\overset{\infty }{\underset{i=0}{\sum }}\overset{j}{\underset{i=0}{\sum }}%
\{[1-p^{(i+2)^{\alpha }}]^{\beta _{2}+\beta _{3}}-[1-p^{(i+1)^{\alpha
}}]^{\beta _{2}+\beta _{3}}\}\{[1-p^{(i+1)^{\alpha }}]^{\beta
_{1}}-[1-p^{i^{\alpha }}]^{\beta _{1}}\} \\
&=&\overset{\infty }{\underset{i=0}{\sum }}\{[1-p^{(i+2)^{\alpha }}]^{\beta
_{2}}-[1-p^{(i+1)^{\alpha }}]^{\beta _{2}}\}[1-p^{(i+1)^{\alpha }}]^{\beta
_{1}+\beta _{3}}.
\end{eqnarray*}

\paragraph{\textit{Result 3:}}

If $X_{1}=\max \{V_{1},V_{3}\}$ and $X_{2}=\max \{V_{2},V_{3}\},$ then $%
X_{1} $ and $X_{2}$ are positive quadrant dependent (PQD).

\paragraph{\textit{Proof:}}

\bigskip \textbf{Case I:} $x_{1}\geq x_{2}$, from (\ref{5}), $P(X_{1}\leq
x_{1},X_{2}\leq x_{2})=F_{1}(x_{1})F_{2}(x_{2})F_{3}(x_{2})$.

Also, from (\ref{1.222}), we get

\begin{equation*}
P(X_{1}\leq x_{1})P(X_{2}\leq
x_{2})=F_{1}(x_{1})F_{3}(x_{1})F_{2}(x_{2})F_{3}(x_{2}).
\end{equation*}%
Therefore, clearly,

\begin{equation*}
P(X_{1}\leq x_{1},X_{2}\leq x_{2})\geq P(X_{1}\leq x_{1})P(X_{2}\leq x_{2}),%
\text{ for all }x_{1}\text{ and }x_{2}.\text{ }
\end{equation*}%
\textbf{Case II:} $x_{1}<x_{2}$. The result can be shown by the same way in
case I.

\paragraph{\textit{Result 4:}}

Suppose $(X_{j1},X_{j2})$ $\sim BDEW(\alpha ,p,\beta _{j1},\beta _{j2},\beta
_{j3})$, for $j=1,...,n$, and they are independently distributed. If $%
Z_{1}=\max \{x_{11},...,x_{n1}\}$ and $Z_{2}=\max \{x_{12},...,x_{n2}\}$,
then

$(Z_{1},Z_{2})$ $\sim BDEW\left( \alpha ,p,\overset{n}{\underset{j=1}{\sum }}%
\beta _{j1},\overset{n}{\underset{j=1}{\sum }}\beta _{j2},\overset{n}{%
\underset{j=1}{\sum }}\beta _{j3}\right) .$

\paragraph{\textit{Proof:}}

The proof can be easily obtained from (\ref{5}).

\section{Reliability Studies}

\subsection{The joint reliability function}

The joint reliability function (RF) of $(X_{1},X_{2})$ can be obtained from
the following relation:%
\begin{equation}
S_{X_{1},X_{2}}(x_{1},x_{2})=1-F_{X_{1}}(x_{1})-F_{X_{2}}(x_{2})+F_{X_{1},X_{2}}(x_{1},x_{2}).
\end{equation}%
The joint RF of $(X_{1},X_{2})$ is given by%
\begin{equation}
S_{X_{1},X_{2}}(x_{1},x_{2})=\left\{ 
\begin{array}{c}
S_{1}(x_{1},x_{2})\text{ \ \ \ \ \ \ \ \ \ \ \ \ \ if \ \ }x_{1}<x_{2} \\ 
S_{2}(x_{1},x_{2})\text{ \ \ \ \ \ \ \ \ \ \ \ \ \ if \ \ }x_{2}<x_{1} \\ 
S_{3}(x)\text{ \ \ \ \ \ \ \ \ \ \ \ \ \ if \ \ }x_{1}=x_{2}=x,%
\end{array}%
\right.   \label{1.27}
\end{equation}%
where

\begin{eqnarray*}
S_{1}(x_{1},x_{2}) &=&1-\text{ }\left[ 1-p^{(x_{2}+1)^{\alpha }}\right]
^{\beta _{2}+\beta _{3}}-\left[ 1-p^{(x_{1}+1)^{\alpha }}\right] ^{\beta
_{1}+\beta _{3}} \\
&&\times \left( 1-\left[ 1-p^{(x_{2}+1)^{\alpha }}\right] ^{\beta
_{2}}\right) ,
\end{eqnarray*}

\begin{eqnarray*}
S_{2}(x_{1},x_{2}) &=&1-\left[ 1-p^{(x_{1}+1)^{\alpha }}\right] ^{\beta
_{1}+\beta _{3}}-\left[ 1-p^{(x_{2}+1)^{\alpha }}\right] ^{\beta _{2}+\beta
_{3}} \\
&&\times \left( 1-\left[ 1-p^{(x_{1}+1)^{\alpha }}\right] ^{\beta
_{1}}\right)
\end{eqnarray*}%
and

\begin{eqnarray*}
S_{3}(x) &=&1-\left[ 1-p^{(x+1)^{\alpha }}\right] ^{\beta _{3}}\times \\
&&\left( \left[ 1-p^{(x+1)^{\alpha }}\right] ^{\beta _{1}}+\left[
1-p^{(x+1)^{\alpha }}\right] ^{\beta _{2}}-\left[ 1-p^{(x+1)^{\alpha }}%
\right] ^{\beta _{1}+\beta _{2}}\right) .
\end{eqnarray*}

\subsection{The joint hazard rate function}

The joint hazard rate function (HRF) of $(X_{1},X_{2})$ can be obtained
using the following relations:%
\begin{equation*}
h_{i}(x_{1},x_{2})=\frac{f_{i}(x_{1},x_{2})}{S_{i}(x_{1},x_{2})},\text{ }%
(i=1,2)\ \text{and }h_{3}(x)=\frac{f_{3}(x)}{S_{3}(x)}.
\end{equation*}%
The joint hazard rate function (HRF) of $(X_{1},X_{2})$ is given by

\begin{equation}
h_{X_{1},X_{2}}(x_{1},x_{2})=\left\{ 
\begin{array}{c}
h_{1}(x_{1},x_{2})\text{ \ \ \ \ \ \ \ \ \ \ \ \ \ if \ \ }x_{1}<x_{2} \\ 
h_{2}(x_{1},x_{2})\text{ \ \ \ \ \ \ \ \ \ \ \ \ \ if \ \ }x_{2}<x_{1} \\ 
h_{3}(x)\text{ \ \ \ \ \ \ \ \ \ \ \ \ \ if \ \ }x_{1}=x_{2}=x,%
\end{array}%
\right.  \label{1.31}
\end{equation}%
where

\begin{equation*}
h_{1}(x_{1},x_{2})=\frac{\left( [1-p^{(x_{1}+1)^{\alpha }}]^{\beta
_{1}+\beta _{3}}-[1-p^{x_{1}^{\alpha }}]^{\beta _{1}+\beta _{3}}\right)
\left( [1-p^{(x_{2}+1)^{\alpha }}]^{\beta _{2}}-[1-p^{x_{2}^{\alpha
}}]^{\beta _{2}}\right) }{1-\text{ }\left[ 1-p^{(x_{2}+1)^{\alpha }}\right]
^{\beta _{2}+\beta _{3}}-\left[ 1-p^{(x_{1}+1)^{\alpha }}\right] ^{\beta
_{1}+\beta _{3}}\left( 1-\left[ 1-p^{(x_{2}+1)^{\alpha }}\right] ^{\beta
_{2}}\right) },
\end{equation*}

\begin{equation*}
h_{2}(x_{1},x_{2})=\frac{\left( [1-p^{(x_{1}+1)^{\alpha }}]^{\beta
_{1}}-[1-p^{x_{1}^{\alpha }}]^{\beta _{1}}\right) \left(
[1-p^{(x_{2}+1)^{\alpha }}]^{\beta _{2}+\beta _{3}}-[1-p^{x_{2}^{\alpha
}}]^{\beta _{2}+\beta _{3}}\right) }{1-\left[ 1-p^{(x_{1}+1)^{\alpha }}%
\right] ^{\beta _{1}+\beta _{3}}-\left[ 1-p^{(x_{2}+1)^{\alpha }}\right]
^{\beta _{2}+\beta _{3}}\left( 1-\left[ 1-p^{(x_{1}+1)^{\alpha }}\right]
^{\beta _{1}}\right) }
\end{equation*}%
and

\begin{equation*}
h_{3}(x)=\frac{[1-p^{(x+1)^{\alpha }}]^{\beta _{1}}\left(
[1-p^{(x+1)^{\alpha }}]^{\beta _{2}+\beta _{3}}-[1-p^{x^{\alpha }}]^{\beta
_{2}+\beta _{3}}\right) -[1-p^{x^{\alpha }}]^{\beta _{1}+\beta _{3}}\left(
[1-p^{(x+1)^{\alpha }}]^{\beta _{2}}-[1-p^{x^{\alpha }}]^{\beta _{2}}\right) 
}{1-\left[ 1-p^{(x+1)^{\alpha }}\right] ^{\beta _{3}}\left( \left[
1-p^{(x+1)^{\alpha }}\right] ^{\beta _{1}}+\left[ 1-p^{(x+1)^{\alpha }}%
\right] ^{\beta _{2}}-\left[ 1-p^{(x+1)^{\alpha }}\right] ^{\beta _{1}+\beta
_{2}}\right) }.
\end{equation*}%
The scatter plot of the hazard rate function of the BDEW distribution is
shown in Figure 2.

\begin{center}
\begin{eqnarray*}
&&\FRAME{itbpFU}{2.5374in}{2.5374in}{0in}{\Qcb{{}}}{}{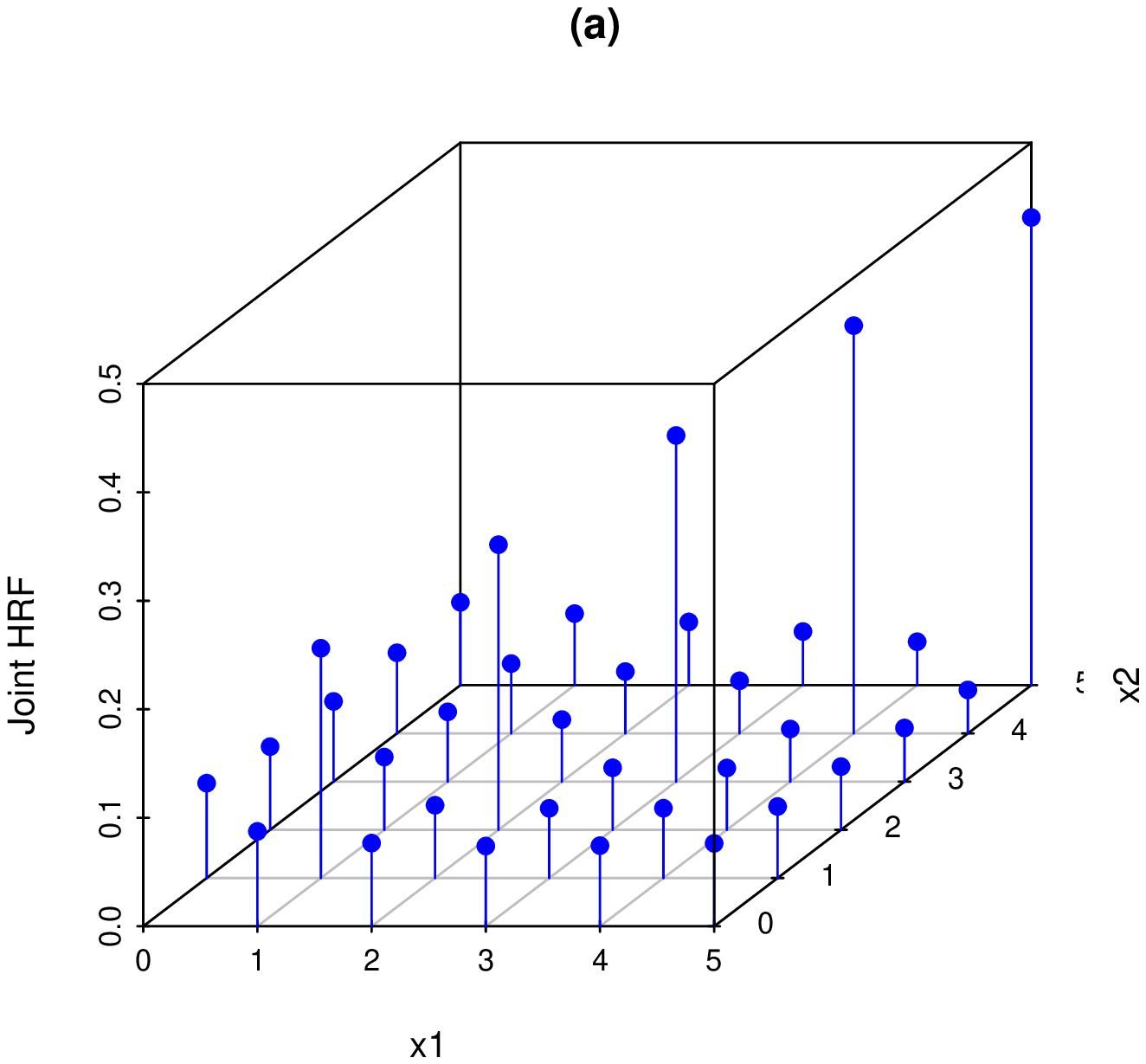}{\special%
{language "Scientific Word";type "GRAPHIC";maintain-aspect-ratio
TRUE;display "USEDEF";valid_file "F";width 2.5374in;height 2.5374in;depth
0in;original-width 5.7519in;original-height 5.7519in;cropleft "0";croptop
"1";cropright "1";cropbottom "0";filename 'h1.eps';file-properties "XNPEU";}}%
\FRAME{itbpFU}{2.5374in}{2.5374in}{0in}{\Qcb{{}}}{}{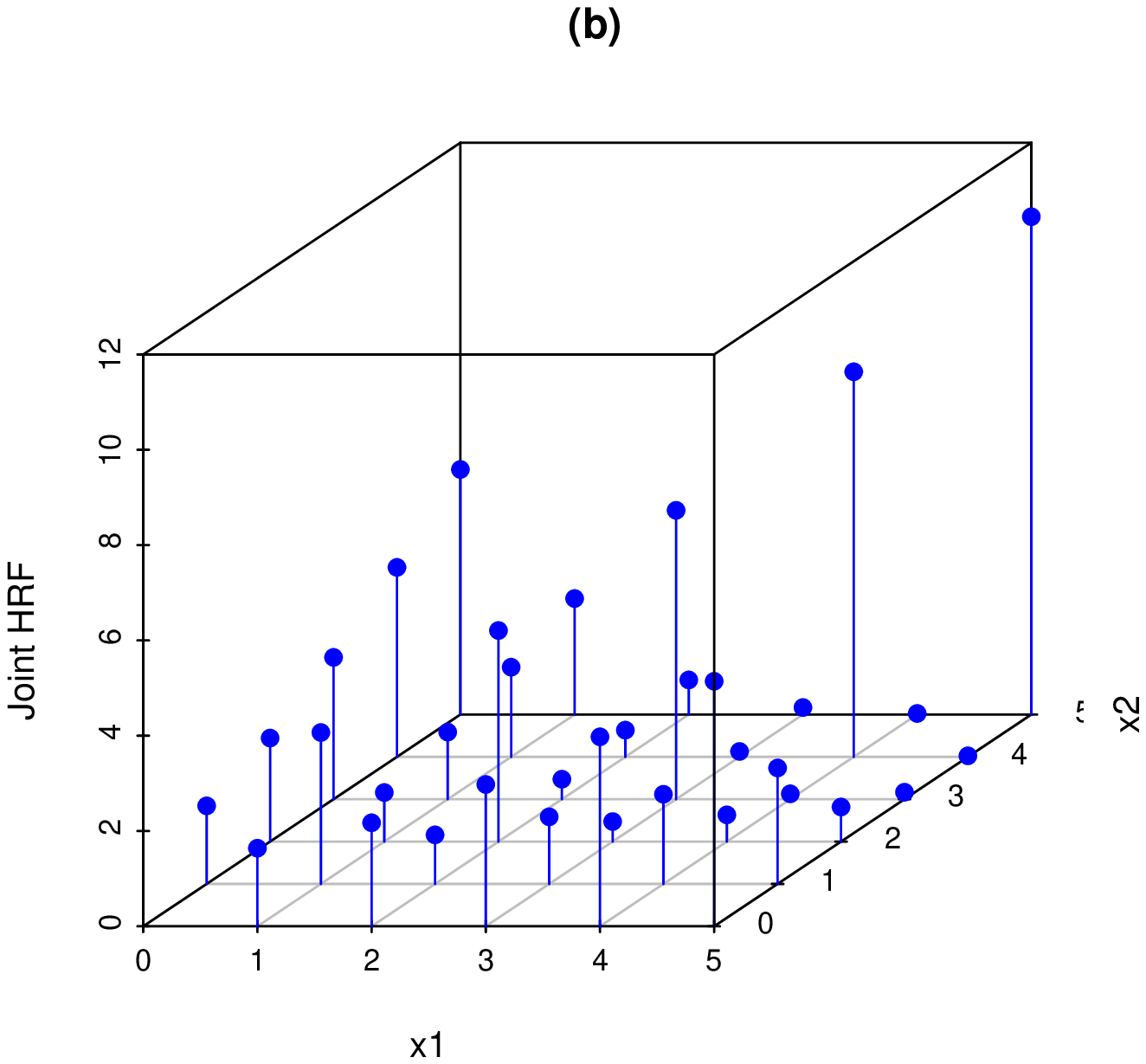}{\special%
{language "Scientific Word";type "GRAPHIC";maintain-aspect-ratio
TRUE;display "USEDEF";valid_file "F";width 2.5374in;height 2.5374in;depth
0in;original-width 5.7519in;original-height 5.7519in;cropleft "0";croptop
"1";cropright "1";cropbottom "0";filename 'h2.eps';file-properties "XNPEU";}}
\\
&&\FRAME{itbpFU}{2.5374in}{2.5374in}{0in}{\Qcb{{}}}{}{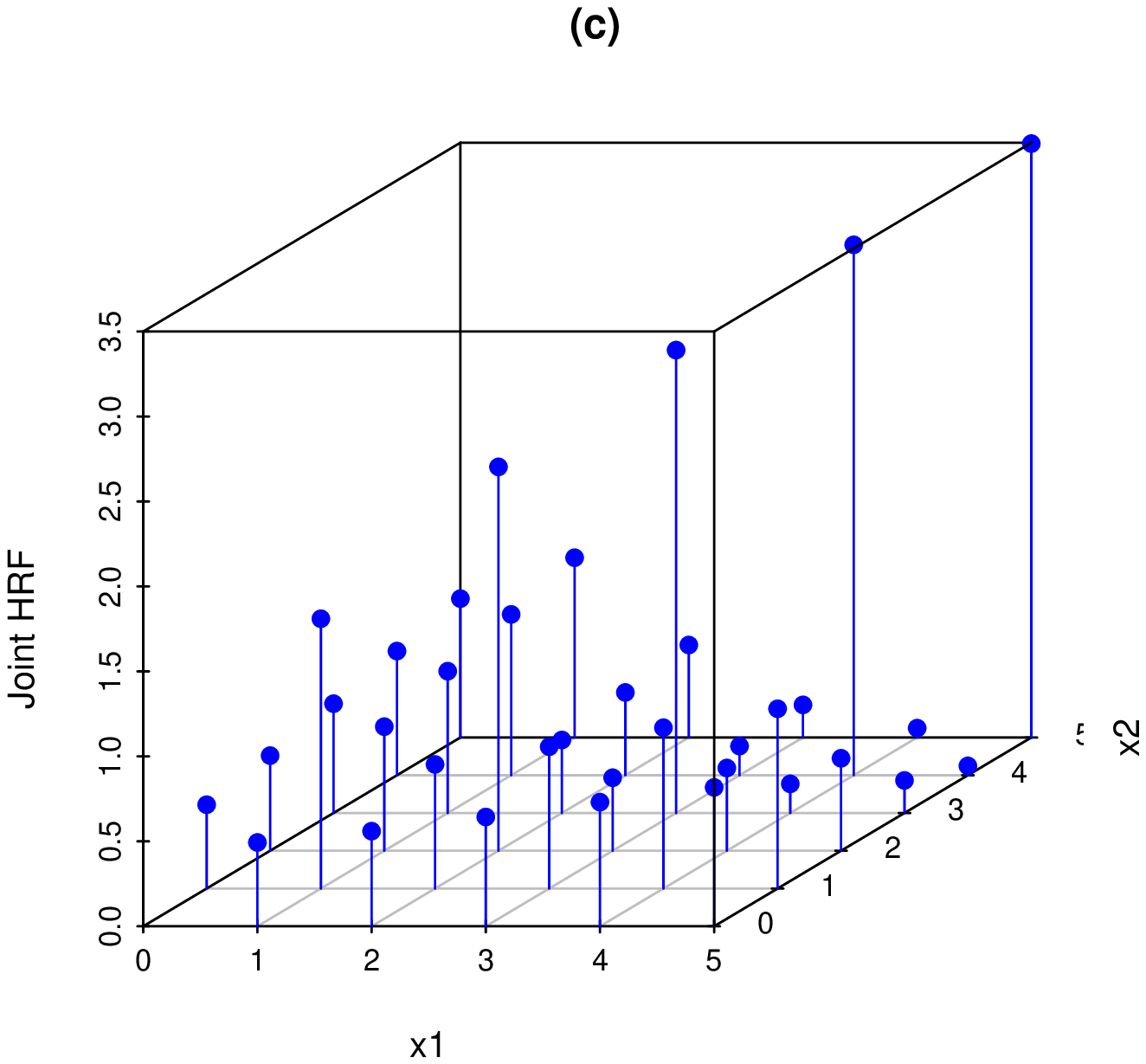}{\special%
{language "Scientific Word";type "GRAPHIC";maintain-aspect-ratio
TRUE;display "USEDEF";valid_file "F";width 2.5374in;height 2.5374in;depth
0in;original-width 5.7519in;original-height 5.7519in;cropleft "0";croptop
"1";cropright "1";cropbottom "0";filename 'h3.eps';file-properties "XNPEU";}}%
\FRAME{itbpFU}{2.5374in}{2.5374in}{0in}{\Qcb{{}}}{}{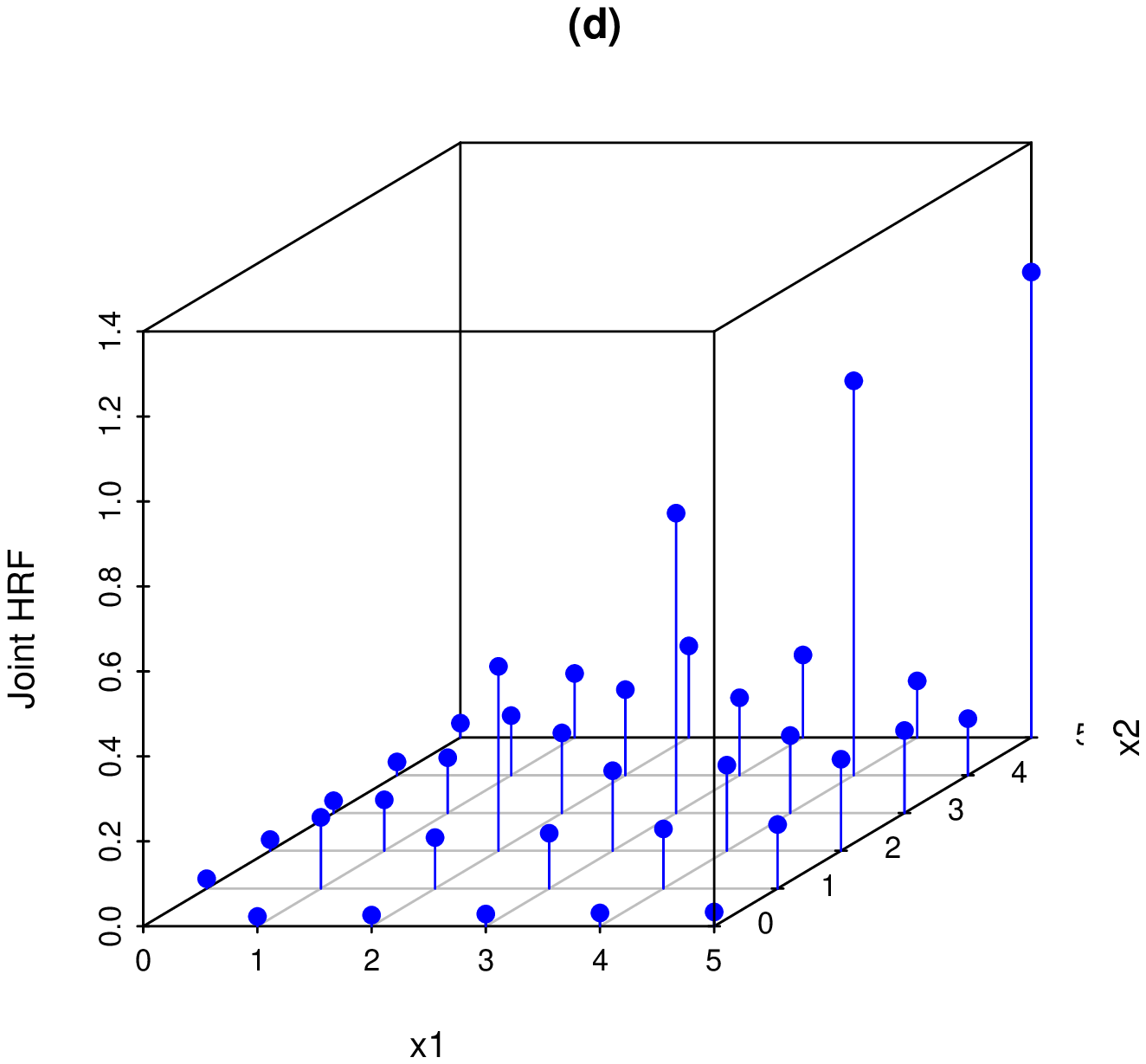}{\special%
{language "Scientific Word";type "GRAPHIC";maintain-aspect-ratio
TRUE;display "USEDEF";valid_file "F";width 2.5374in;height 2.5374in;depth
0in;original-width 5.7519in;original-height 5.7519in;cropleft "0";croptop
"1";cropright "1";cropbottom "0";filename 'h4.eps';file-properties "XNPEU";}}
\end{eqnarray*}%
\textbf{Figure 2.} Scatter plot of the joint HRF of the BDEW distribution
for different values of its parameter vector $\mathbf{\Omega }=(\alpha
,p,\beta _{1},\beta _{2},\beta _{3})$ : \textbf{(a)} $\mathbf{\Omega }%
=(1.5,0.9,0.3,0.3,0.3),$ \textbf{(b)} $\mathbf{\Omega }%
=(1.5,0.5,0.3,0.3,0.3),$ \textbf{(c)} $\mathbf{\Omega }=(1.3,0.5,0.7,0.7,0.7)
$ and \textbf{(d)} $\mathbf{\Omega }=(1.3,0.7,1.2,1.2,1.2)$.
\end{center}

\section{Maximum Likelihood Estimation}

In this section, we use the method of maximum likelihood to estimate the
unknown parameters $\alpha ,p,\beta _{1},\beta _{2}$ and $\beta _{3}$ of
the\ BDEW distribution. Suppose that, we have a sample of size n, of the
form $\left\{ (x_{11},x_{21}),(x_{12},x_{22}),...,(x_{1n},x_{2n})\right\} $
from the BDEW distribution. We use the following notations: $%
I_{1}=\{x_{1j}<x_{2j}\},$ $I_{2}=\{x_{2j}<x_{1j}\},$ $I_{3}=%
\{x_{1j}=x_{2j}=x_{j}\},$ $I=I_{1}\cup I_{2}\cup I_{3},$ $\left\vert
I_{1}\right\vert =n_{1},$ $\left\vert I_{2}\right\vert =n_{2},$ $\left\vert
I_{3}\right\vert =n_{3}$ and $n=n_{1}+n_{2}+n_{3}.$

Based on the observations, the likelihood function is given by

\begin{equation*}
l(\mathbf{\Omega })=\underset{j=1}{\overset{n_{1}}{\prod }}%
f_{1}(x_{1j},x_{2j})\underset{j=1}{\overset{n_{2}}{\prod }}%
f_{2}(x_{1j},x_{2j})\underset{j=1}{\overset{n_{3}}{\prod }}f_{3}(x_{j}).
\end{equation*}%
The log-likelihood function becomes

\begin{eqnarray}
L(\mathbf{\Omega }) &=&\overset{n_{1}}{\underset{j=1}{\sum }}\ln \left(
g_{1}(x_{1j};\beta _{1}+\beta _{3})\right) +\overset{n_{1}}{\underset{j=1}{%
\sum }}\ln \left( g_{1}(x_{2j};\beta _{2})\right)   \notag \\
&&+\overset{n_{2}}{\underset{j=1}{\sum }}\ln \left( g_{1}(x_{1j};\beta
_{1})\right) +\overset{n_{2}}{\underset{j=1}{\sum }}\ln \left(
g_{1}(x_{2j};\beta _{2}+\beta _{3})\right)   \notag \\
&&+\overset{n_{3}}{\underset{j=1}{\sum }}\ln \left( [1-p^{(x_{j}+1)^{\alpha
}}]^{\beta _{1}}g_{1}(x_{j};\beta _{2}+\beta _{3})-[1-p^{x_{j}{}^{\alpha
}}]^{\beta _{1}+\beta _{3}}g_{1}(x_{j}+1;\beta _{2})\right) ,  \notag \\
&&  \label{1.36}
\end{eqnarray}%
where%
\begin{equation*}
g_{1}(x;\beta )=[1-p^{(x+1)^{\alpha }}]^{\beta }-[1-p^{x{}^{\alpha
}}]^{\beta }.
\end{equation*}

The MLEs of the parameters $\alpha ,p,\beta _{1},\beta _{2}$ and $\beta _{3}$
can be obtained by computing the first partial derivatives of (\ref{1.36})
with respect to $\alpha ,p,\beta _{1},\beta _{2}$ and $\beta _{3}$, putting
the results equal zeros. We get the likelihood equations as in the following
form%
\begin{eqnarray}
\frac{\partial L}{\partial \alpha } &=&\overset{n_{1}}{\underset{j=1}{\sum }}%
\frac{g_{4}(x_{1j}+1;\beta _{1}+\beta _{3})-g_{4}(x_{1j};\beta _{1}+\beta
_{3})}{g_{1}(x_{1j};\beta _{1}+\beta _{3})}+\overset{n_{1}}{\underset{j=1}{%
\sum }}\frac{g_{4}(x_{2j}+1;\beta _{2})-g_{4}(x_{2j};\beta _{2})}{%
g_{1}(x_{2j};\beta _{2})}  \notag \\
&&+\overset{n_{2}}{\underset{j=1}{\sum }}\frac{g_{4}(x_{1j}+1;\beta
_{1})-g_{4}(x_{1j};\beta _{1})}{g_{1}(x_{1j};\beta _{1})}+\overset{n_{2}}{%
\underset{j=1}{\sum }}\frac{g_{4}(x_{2j}+1;\beta _{2}+\beta
_{3})-g_{4}(x_{2j};\beta _{2}+\beta _{3})}{g_{1}(x_{2j};\beta _{2}+\beta
_{3})}  \notag \\
&&+\overset{n_{3}}{\underset{j=1}{\sum }}\frac{[1-p^{(x_{j}+1)^{\alpha
}}]^{\beta _{1}}\left( g_{4}(x_{j}+1;\beta _{2}+\beta
_{3})-g_{4}(x_{j};\beta _{2}+\beta _{3})\right) +g_{4}(x_{j}+1;\beta
_{1})g_{1}(x_{j};\beta _{2}+\beta _{3})}{[1-p^{(x_{j}+1)^{\alpha }}]^{\beta
_{1}}g_{1}(x_{j};\beta _{2}+\beta _{3})-[1-p^{x_{j}{}^{\alpha }}]^{\beta
_{1}+\beta _{3}}g_{1}(x_{j}+1;\beta _{2})}  \notag \\
&&-\overset{n_{3}}{\underset{j=1}{\sum }}\frac{[1-p^{x_{j}{}^{\alpha
}}]^{\beta _{1}+\beta _{3}}\left( g_{4}(x_{j}+1;\beta
_{2})-g_{4}(x_{j};\beta _{2})\right) +g_{4}(x_{j};\beta _{1}+\beta
_{3})g_{1}(x_{j};\beta _{2})}{[1-p^{(x_{j}+1)^{\alpha }}]^{\beta
_{1}}g_{1}(x_{j};\beta _{2}+\beta _{3})-[1-p^{x_{j}{}^{\alpha }}]^{\beta
_{1}+\beta _{3}}g_{1}(x_{j}+1;\beta _{2})},  \label{1.37}
\end{eqnarray}%
\begin{eqnarray}
\frac{\partial L}{\partial p} &=&\overset{n_{1}}{\underset{j=1}{\sum }}\frac{%
g_{3}(x_{1j}+1;\beta _{1}+\beta _{3})-g_{3}(x_{1j};\beta _{1}+\beta _{3})}{%
g_{1}(x_{1j};\beta _{1}+\beta _{3})}+\overset{n_{1}}{\underset{j=1}{\sum }}%
\frac{g_{3}(x_{2j}+1;\beta _{2})-g_{3}(x_{2j};\beta _{2})}{%
g_{1}(x_{2j};\beta _{2})}  \notag \\
&&+\overset{n_{2}}{\underset{j=1}{\sum }}\frac{g_{3}(x_{1j}+1;\beta
_{1})-g_{3}(x_{1j};\beta _{1})}{g_{1}(x_{1j};\beta _{1})}+\overset{n_{2}}{%
\underset{j=1}{\sum }}\frac{g_{3}(x_{2j}+1;\beta _{2}+\beta
_{3})-g_{3}(x_{2j};\beta _{2}+\beta _{3})}{g_{1}(x_{2j};\beta _{2}+\beta
_{3})}  \notag \\
&&+\overset{n_{3}}{\underset{j=1}{\sum }}\frac{[1-p^{(x_{j}+1)^{\alpha
}}]^{\beta _{1}}\left( g_{3}(x_{j}+1;\beta _{2}+\beta
_{3})-g_{3}(x_{j};\beta _{2}+\beta _{3})\right) +g_{3}(x_{j}+1;\beta
_{1})g_{1}(x_{j};\beta _{2}+\beta _{3})}{[1-p^{(x_{j}+1)^{\alpha }}]^{\beta
_{1}}g_{1}(x_{j},;\beta _{2}+\beta _{3})-[1-p^{x_{j}{}^{\alpha }}]^{\beta
_{1}+\beta _{3}}g_{1}(x_{j}+1;\beta _{2})}  \notag \\
&&-\overset{n_{3}}{\underset{j=1}{\sum }}\frac{[1-p^{x_{j}{}^{\alpha
}}]^{\beta _{1}+\beta _{3}}\left( g_{3}(x_{j}+1;\beta
_{2})-g_{3}(x_{j};\beta _{2})\right) +g_{3}(x_{j};\beta _{1}+\beta
_{3})g_{1}(x_{j};\beta _{2})}{[1-p^{(x_{j}+1)^{\alpha }}]^{\beta
_{1}}g_{1}(x_{j};\beta _{2}+\beta _{3})-[1-p^{x_{j}{}^{\alpha }}]^{\beta
_{1}+\beta _{3}}g_{1}(x_{j}+1;\beta _{2})},  \label{1.38}
\end{eqnarray}%
\begin{eqnarray}
\frac{\partial L}{\partial \beta _{1}} &=&\overset{n_{1}}{\underset{j=1}{%
\sum }}\frac{g_{2}(x_{1j}+1;\beta _{1}+\beta _{3})-g_{2}(x_{1j};\beta
_{1}+\beta _{3})}{g_{1}(x_{1j};\beta _{1}+\beta _{3})}+\overset{n_{2}}{%
\underset{j=1}{\sum }}\frac{g_{2}(x_{1j}+1;\beta _{1})-g_{2}(x_{1j},\beta
_{1})}{g_{1}(x_{1j};\beta _{1})}  \notag \\
&&+\overset{n_{3}}{\underset{j=1}{\sum }}\frac{g_{2}(x_{j}+1;\beta
_{1})g_{1}(x_{j};\beta _{2}+\beta _{3})-g_{2}(x_{j};\beta _{1}+\beta
_{3})g_{1}(x_{j};\beta _{2})}{[1-p^{(x_{j}+1)^{\alpha }}]^{\beta
_{1}}g_{1}(x_{j};\beta _{2}+\beta _{3})-[1-p^{x_{j}{}^{\alpha }}]^{\beta
_{1}+\beta _{3}}g_{1}(x_{j}+1;\beta _{2})},  \label{1.39}
\end{eqnarray}%
\begin{eqnarray}
\frac{\partial L}{\partial \beta _{2}} &=&\overset{n_{1}}{\underset{j=1}{%
\sum }}\frac{g_{2}(x_{2j}+1;\beta _{2})-g_{2}(x_{2j};\beta _{1})}{%
g_{1}(x_{2j};\beta _{2})}+\overset{n_{2}}{\underset{j=1}{\sum }}\frac{%
g_{2}(x_{2j}+1;\beta _{2}+\beta _{3})-g_{2}(x_{2j};\beta _{2}+\beta _{3})}{%
g_{1}(x_{2j};\beta _{2}+\beta _{3})}  \notag \\
&&+\overset{n_{3}}{\underset{j=1}{\sum }}\frac{[1-p^{(x_{j}+1)^{\alpha
}}]^{\beta _{1}}\left( g_{2}(x_{j}+1;\beta _{2}+\beta
_{3})-g_{2}(x_{j};\beta _{2}+\beta _{3})\right) }{[1-p^{(x_{j}+1)^{\alpha
}}]^{\beta _{1}}g_{1}(x_{j};\beta _{2}+\beta _{3})-[1-p^{x_{j}{}^{\alpha
}}]^{\beta _{1}+\beta _{3}}g_{1}(x_{j}+1;\beta _{2})}  \notag \\
&&-\overset{n_{3}}{\underset{j=1}{\sum }}\frac{[1-p^{x_{j}{}^{\alpha
}}]^{\beta _{1}+\beta _{3}}\left( g_{2}(x_{j}+1;\beta
_{2})-g_{2}(x_{j};\beta _{2})\right) }{[1-p^{(x_{j}+1)^{\alpha }}]^{\beta
_{1}}g_{1}(x_{j};\beta _{2}+\beta _{3})-[1-p^{x_{j}{}^{\alpha }}]^{\beta
_{1}+\beta _{3}}g_{1}(x_{j}+1;\beta _{2})},  \label{1.40}
\end{eqnarray}%
and%
\begin{eqnarray}
\frac{\partial L}{\partial \beta _{3}} &=&\overset{n_{1}}{\underset{j=1}{%
\sum }}\frac{g_{2}(x_{1j}+1;\beta _{1}+\beta _{3})-g_{2}(x_{1j};\beta
_{1}+\beta _{3})}{g_{1}(x_{1j};\beta _{1}+\beta _{3})}+\overset{n_{2}}{%
\underset{j=1}{\sum }}\frac{g_{2}(x_{2j}+1;\beta _{2}+\beta
_{3})-g_{2}(x_{2j};\beta _{2}+\beta _{3})}{g_{1}(x_{2j};\beta _{2}+\beta
_{3})}  \notag \\
&&+\overset{n_{3}}{\underset{j=1}{\sum }}\frac{[1-p^{(x_{j}+1)^{\alpha
}}]^{\beta _{1}}\left( g_{2}(x_{j}+1;\beta _{2}+\beta
_{3})-g_{2}(x_{j};\beta _{2}+\beta _{3})\right) -g_{2}(x_{j};\beta
_{1}+\beta _{3})g_{1}(x_{j};\beta _{2})}{[1-p^{(x_{j}+1)^{\alpha }}]^{\beta
_{1}}g_{1}(x_{j};\beta _{2}+\beta _{3})-[1-p^{x_{j}{}^{\alpha }}]^{\beta
_{1}+\beta _{3}}g_{1}(x_{j}+1;\beta _{2})},  \notag \\
&&  \label{1.41}
\end{eqnarray}%
where%
\begin{eqnarray*}
g_{2}(x;\beta ) &=&[1-p^{x{}^{\alpha }}]^{\beta }\ln (1-p^{x{}^{\alpha }}),
\\
g_{3}(x;\beta ) &=&-\beta x^{\alpha }p^{x^{\alpha }-1}[1-p^{x{}^{\alpha
}}]^{\beta -1}, \\
g_{4}(x;\beta ) &=&-\beta \ln (x)x^{\alpha }p^{x^{\alpha }}\ln
(p)[1-p^{x{}^{\alpha }}]^{\beta -1}.
\end{eqnarray*}%
The MLEs of the parameters $\alpha ,$ $p,$ $\beta _{1},$ $\beta _{2}$ and $%
\beta _{3}$ can be obtained by solving the above system of five non-linear
equations from (\ref{1.37}) to (\ref{1.41}). The solution of thes equations
are not easy to solve, so we need a numerical technique to get the MLEs.

\section{Data Analysis}

In this section, we explain the experimental importance of the BDEW
distribution using two applications to real data sets. In each data, we
shall compare the fits of the BDEW distribution with some competitive models
such as BGDR, BDGE, bivariate geometric (BG) distribution which introduced
by Basu and Dhar (1995), Independent bivariate Poisson (IBP), bivariate
Poisson (BP) by Holgate (1964) and bivariate Poisson with minimum operator
(BP$_{\min }$) by Lee and Cha (2015). The tested distributions are compared
using some criteria namely, the maximized Log-Likelihood (\ $-L$\ ), Akaike
information criterion (AIC), bayesian information Criterion (BIC), corrected
Akaike information criterion (CAIC) and Hannan-Quinn information criterion
(HQIC). Further, we can use the Pearson's chi-square goodness-of-fit test
(see, Lawless(1982)) for grouped data to test the goodness of fit of a
proposed bivariate distribution. But the sample size must be sufficiently
large in order to apply this test. For this reason, we did not use this test
in the two data sets analyzed here.

\paragraph{Data 1:}

The data set given in Table 1 consists of a football match score in Italian
football match (Serie A) during 1996 to 2011, between ACF Fiorentina($X_{1}$%
) and Juventus($X_{2}$). For previous studies for this data, see Lee and Cha
(2015), Nekoukhou and Kundu (2017) and Xiao et al. (2017). Unfortunately,
there is an error in the fourth and tenth observations of the variable $X_{1}
$ in Nekoukhou and Kundu (2017) paper, which raised the results of this
study. The data source is (
http://www.worldfootball.net/competition/ita-serie-a/). The MLEs, -L, AIC,
CAIC, and HQIC values for\ IBP, BP, BP$_{\min }$, BG and BDEW distributions
are obtained in Table 2.%
\begin{eqnarray*}
&&\text{\textbf{Table 1. }The\textbf{\ }score data between ACF Fiorentina(X}%
_{1}\text{) and Juventus(X}_{2}\text{).} \\
&&%
\begin{tabular}{|c|c|c|c|c|c|c|c|}
\hline\hline
\textbf{Obs.} & \ \ \ \textbf{Match Date \ \ } & $X_{1}$ & $X_{2}$ & \ \ \ 
\textbf{Obs.} & \ \ \ \textbf{Match Date \ \ } & $X_{1}$ & $X_{2}$ \\ 
\hline\hline
\textbf{1} & 25/10/2011 & $1$ & $2$ & \textbf{14} & 16/02/2002 & $1$ & $2$
\\ \hline
\textbf{2} & 17/04/2011 & $0$ & $0$ & \textbf{15} & 19/12/2001 & $1$ & $1$
\\ \hline
\textbf{3} & 27/11/2010 & $1$ & $1$ & \textbf{16} & 12/05/2001 & $1$ & $3$
\\ \hline
\textbf{4} & 06/03/2010 & $1$ & $2$ & \textbf{17} & 06/01/2001 & $3$ & $3$
\\ \hline
\textbf{5} & 17/10/2009 & $1$ & $1$ & \textbf{18} & 21/04/2000 & $0$ & $1$
\\ \hline
\textbf{6} & 24/01/2009 & $0$ & $1$ & \textbf{19} & 18/12/1999 & $1$ & $1$
\\ \hline
\textbf{7} & 31/08/2008 & $1$ & $1$ & \textbf{20} & 24/04/1999 & $1$ & $2$
\\ \hline
\textbf{8} & 02/03/2008 & $3$ & $2$ & \textbf{21} & 12/12/1998 & $1$ & $0$
\\ \hline
\textbf{9} & 07/10/2007 & $1$ & $1$ & \textbf{22} & 21/02/1998 & $3$ & $0$
\\ \hline
\textbf{10} & 09/04/2006 & $1$ & $1$ & \textbf{23} & 04/10/1997 & $1$ & $2$
\\ \hline
\textbf{11} & 04/12/2005 & $1$ & $2$ & \textbf{24} & 22/02/1997 & $1$ & $1$
\\ \hline
\textbf{12} & 09/04/2005 & $3$ & $3$ & \textbf{25} & 28/09/1996 & $0$ & $1$
\\ \hline
\textbf{13} & 10/11/2004 & $0$ & $1$ & \textbf{26} & 23/03/1996 & $0$ & $1$
\\ \hline\hline
\end{tabular}%
\end{eqnarray*}

\begin{eqnarray*}
&&\text{\textbf{Table 2. } The MLEs, -L, AIC, CAIC, BIC and HQIC values for
data set 1}. \\
&&%
\begin{tabular}{c|c|c|c|c|c}
\hline\hline
& \multicolumn{5}{|c}{\textbf{Distributions}} \\ \hline
\textbf{Statistic} & \textbf{\ \ \ \ IBP \ \ \ \ } & \textbf{\ \ \ \ BP \ \
\ \ } & \ \ \textbf{BP}$_{\min }$ \ \  & \textbf{\ \ \ \ \ BG \ \ \ \ \ } & 
\ \ \textbf{BDEW \ \ } \\ \hline
$\overset{\wedge }{\alpha }$ & $--$ & $--$ & $--$ & $0.343$ & $0.922$ \\ 
\hline
$\overset{\wedge }{p}$ & $--$ & $--$ & $--$ & $0.553$ & $0.172$ \\ \hline
$\mathbf{\ }\widehat{\beta }_{1}$ & $1.08$ & $1.08$ & $1.36$ & $--$ & $4.315$
\\ \hline
$\mathbf{\ }\widehat{\beta }_{2}$ & $1.38$ & $1.38$ & $2.10$ & $--$ & $9.656$
\\ \hline
$\mathbf{\ }\widehat{\beta }_{3}$ & $--$ & $0.70$ & $2.27$ & $--$ & $2.892$
\\ \hline
\textbf{-L} & $67.60$ & $64.92$ & $64.22$ & $76.36$ & $60.89$ \\ \hline
\textbf{AIC} & $139.21$ & $135.83$ & $134.44$ & $156.72$ & $131.7$ \\ \hline
\textbf{CAIC} & $139.72$ & $136.93$ & $135.53$ & $157.24$ & $134.83$ \\ 
\hline
\textbf{BIC} & $141.72$ & $139.61$ & $138.21$ & $159.24$ & $138.09$ \\ \hline
\textbf{HQIC} & $139.92$ & $136.93$ & $135.53$ & $157.44$ & $133.61$ \\ 
\hline\hline
\end{tabular}%
\end{eqnarray*}%
\ We can conclude that, the BDEW distribution provides a better fit than the
other tested distributions, because the BDEW distribution has the smallest
values for -L, AIC, CAIC, BIC and HQIC among all tested distributions.

\paragraph{Data 2:}

The data set given in Table 3 is taken from a video recording in 1995 IX
World Cup diving Championship, Atlanta, Georgia by NBC sports TV. This data
consists of a\ scores given by seven judges from seven varies countries
recorded in a video that starts at the end of the fourth round, which is a
random start, taken from NBC sports TV. The score given by each judge is
taken positive integer values. In this data, the random variable $X_{1}$
represent the max score between Asian and Caucasus countries, but $X_{2}$
represent the max score between western countries, United Kingdom, Canada,
Australia, Iceland and France. The score corresponding to the dive of
Michael Murphy (Obs. number 3) was not displayed by NBC sports TV. For
previous studies for this data, see Jing and Dhar (2012). In Table 4, the
MLEs, -L, AIC, CAIC, and HQIC values for\ BG, BDGE, BDERe and BDEW
distributions are obtained.%
\begin{eqnarray*}
&&\text{\textbf{Table 3. }The\textbf{\ }video recording data by NBC sports
TV.} \\
&&%
\begin{tabular}{|c|c|c|c||c|c|c|c|}
\hline\hline
{\small Obs.} & \textbf{Diver} & ${\small X}_{1}$ & ${\small X}_{2}$ & 
{\small Obs.} & \textbf{Diver} & ${\small X}_{1}$ & ${\small X}_{2}$ \\ 
\hline
{\small 1} & \multicolumn{1}{|l|}{\small Sun Shuwei, China} & ${\small 19}$
& ${\small 19}$ & {\small 11} & \multicolumn{1}{|l|}{\small Sun Shuwei, China%
} & ${\small 15}$ & ${\small 16}$ \\ \hline
{\small 2} & \multicolumn{1}{|l|}{\small David Pichler, USA} & ${\small 15}$
& ${\small 15}$ & {\small 12} & \multicolumn{1}{|l|}{\small David Pichler,
USA} & ${\small 15}$ & ${\small 15}$ \\ \hline
{\small 3} & \multicolumn{1}{|l|}{\small Michael Murphy, Australia} & $%
{\small --}$ & ${\small --}$ & {\small 13} & \multicolumn{1}{|l|}{\small Jan
Hempel, Germany} & ${\small 17}$ & ${\small 18}$ \\ \hline
{\small 4} & \multicolumn{1}{|l|}{\small Jan Hempel, Germany} & ${\small 13}$
& ${\small 14}$ & {\small 14} & \multicolumn{1}{|l|}{\small Roman Volodkuv,
Ukraine} & ${\small 16}$ & ${\small 16}$ \\ \hline
{\small 5} & \multicolumn{1}{|l|}{\small Roman Volodkuv, Ukraine} & ${\small %
11}$ & ${\small 12}$ & {\small 15} & \multicolumn{1}{|l|}{\small Sergei
Kudrevich, Belarus} & ${\small 12}$ & ${\small 13}$ \\ \hline
{\small 6} & \multicolumn{1}{|l|}{\small Sergei Kudrevich, Belarus} & $%
{\small 14}$ & ${\small 14}$ & {\small 16} & \multicolumn{1}{|l|}{\small %
Patrick Je rey, USA} & ${\small 14}$ & ${\small 14}$ \\ \hline
{\small 7} & \multicolumn{1}{|l|}{\small Patrick Je rey, USA} & ${\small 15}$
& ${\small 14}$ & {\small 17} & \multicolumn{1}{|l|}{\small Valdimir
Timoshinin, Russia} & ${\small 12}$ & ${\small 13}$ \\ \hline
{\small 8} & \multicolumn{1}{|l|}{\small Valdimir Timoshinin, Russia} & $%
{\small 13}$ & ${\small 16}$ & {\small 18} & \multicolumn{1}{|l|}{\small %
Dimitry Sautin, Russia} & ${\small 17}$ & ${\small 18}$ \\ \hline
{\small 9} & \multicolumn{1}{|l|}{\small Dimitry Sautin, Russia} & ${\small 7%
}$ & ${\small 5}$ & {\small 19} & \multicolumn{1}{|l|}{\small Xiao Hailiang,
China} & ${\small 9}$ & ${\small 10}$ \\ \hline
{\small 10} & \multicolumn{1}{|l|}{\small Xiao Hailiang, China} & ${\small 13%
}$ & ${\small 13}$ & {\small 20} & \multicolumn{1}{|l|}{\small Sun Shuwei,
China} & ${\small 18}$ & ${\small 18}$ \\ \hline
\end{tabular}%
\end{eqnarray*}

\begin{eqnarray*}
&&\text{\textbf{Table 4. }The MLEs, -L, AIC, CAIC, BIC and HQIC values for
data set 2.} \\
&&%
\begin{tabular}{c|c|c|c|c}
\hline\hline
& \multicolumn{4}{|c}{\textbf{Distributions}} \\ \hline
\textbf{Statistic} & \textbf{\ \ \ \ \ \ BG\ \ \ \ \ \ } & \ \ \ \ \textbf{\
B\textbf{D}GE \ \ \ \ } & \textbf{\ \ \ \ BDER \ \ \ \ \ } & \ \ \ \textbf{%
BDEW \ \ } \\ \hline
$\overset{\wedge }{\alpha }$ & $0.5277$ & $--$ & $--$ & $3.5239$ \\ \hline
$\overset{\wedge }{p}$ & $0.9357$ & $0.7613$ & $0.9893$ & $0.9998$ \\ \hline
$\mathbf{\ }\widehat{\beta }_{1}$ & $--$ & $8.6046$ & $1.4968$ & $0.5327$ \\ 
\hline
$\mathbf{\ }\widehat{\beta }_{2}$ & $--$ & $17.981$ & $3.3389$ & $1.2491$ \\ 
\hline
$\mathbf{\ }\widehat{\beta }_{3}$ & $--$ & $24.116$ & $4.3561$ & $1.5922$ \\ 
\hline
\textbf{-L} & $121.584$ & $90.959$ & $86.737$ & $84.056$ \\ \hline
\textbf{AIC} & $247.168$ & $189.917$ & $181.474$ & $178.111$ \\ \hline
\textbf{CAIC} & $247.918$ & $192.774$ & $184.331$ & $182.726$ \\ \hline
\textbf{BIC} & $249.057$ & $193.695$ & $185.251$ & $182.834$ \\ \hline
\textbf{HQIC} & $247.488$ & $190.556$ & $182.113$ & $178.911$ \\ \hline\hline
\end{tabular}%
\end{eqnarray*}%
It is clear that, the BDEW distribution has the smallest values for -L, AIC,
CAIC and HQIC among all tested distributions. So, the BDEW distribution
provides a better fit than the other tested distributions.

\section{Conclusions}

\ \ \ \ In this paper we introduced a new five parameters bivariate discrete
distribution called the bivariate discrete exponentiated Weibull
distribution. Some statistical and reliability properties of the proposed
discrete model have been derived. The maximum likelihood estimators of the
parameters are deduced. Two real data sets have been analyzed using the new
discrete distribution compared with another famous distributions. It is
clear from the comparison that the new distribution is the best distribution
for fitting the data sets from among the all tested distributions.

\bigskip \setlength{\parindent}{0in}


\begin{thebibliography}{99}
\bibitem{16} Alamatsaz, M. H., Dey, S., Dey, T., and Harandi, S. Shams.
(2016). Discrete generalized Rayleigh distribution, Pakistan journal of
statistics, 32(1), 1-20.

\bibitem{32} Basu, A. P. and Dhar, S. K. (1995). Bivariate geometric
distribution. Journal of applied statistical science, 2, 33-34.

\bibitem{4} Bebbington, M., Lai, C. D., and Zitikis, R. (2007). A flexible
Weibull extension. Reliability engineering and system safety, 92, 719-726.

\bibitem{23} El-Bassiouny, A. H., EL-Damcese, M., Abdelfattah, M., and
Eliwa, M. S. (2016). Bivariate exponentaited generalized Weibull-Gompertz
distribution. Journal of applied probability and statistics, 11(1), 25-46.

\bibitem{39} El-Bassiouny, A. H., Medhat EL-Damcese, Abdelfattah Mustafa,
and Eliwa, M. S. (2017). Exponentiated generalized Weibull-Gompertz
distribution with application in survival analysis. Journal of statistics
applications and probability, 6(1), 7-16.

\bibitem{5} El-Gohary, A., EL-Bassiouny, A. H., and El-Morshedy, M. (2015).
Exponentiated flexible Weibull extension distribution. International journal
of mathematics and its applications, 3(A ), 1-12.

\bibitem{25} El-Gohary, A., El-Bassiouny, A. H., and El-Morshedy, M.,
(2016). Bivariate exponentiated modified Weibull extension distribution.
Journal of statistics applications and probability, 5(1), 67-78.

\bibitem{7} El-Morshedy, M., El-Bassiouny, A. H., and El-Gohary, A., (2017).
Exponentiated inverse flexible Weibull extension distribution. Journal of
statistics applications and probability, 6(1), 169-183.

\bibitem{21} El-Sherpieny, E. A., Ibrahim, S. A., and Bedar, R. E. (2013). A
new bivariate generalized Gompertz distribution. Asian journal of applied
sciences, 1(4): 141-150.

\bibitem{13} Gomez-Deniz, E., and Calderin-Ojeda, E. (2011). The discrete
Lindley distribution: properties and applications. Journal of statistical
computation and simulation, 81(11), 1405-1416.
doi:10.1080/00949655.2010.487825.

\bibitem{34} Holgate, B. (1964). Estimation for the bivariate Poisson
distribution,\ Biometrika, 51, 241-245.

\bibitem{38} Jing, Li, and Sunil, K. Dhar. (2012). Modeling with bivariate
Geometric distributions, Communications in Statistics - Theory and Methods,
42(2), 252-266, DOI: 10.1080/03610926.2011.579704.

\bibitem{18} Jose, K.K., Ristic, M.M., and Ancy, J. (2009). Marshall- Olkin
bivariate Weibull distributions and processes. Statistical papers, DOI
10.1007/s00362-009-0287-8.

\bibitem{29} Kemp, A. W. (2013). New discrete appell and humbert
distributions with relevance to bivariate accident data. Journal of
multivariate analysis, 113, 2-6.

\bibitem{27} Kocherlakota, S., and Kocherlakota, K. (1992). Bivariate
discrete distributions. Marcel dekker, New York.

\bibitem{12} Krishna, H., and Pundir, P. S. (2009). Discrete Burr and
discrete Pareto distributions. Statistical methodology, 6(2), 177-188.
doi:10.1016/j.stamet.2008.07.001.

\bibitem{28} Kumar, C. S. (2008). A unified approach to bivariate discrete
distributions. Metrika, 67, 113-123.

\bibitem{19} Kundu, D., and Gupta, R. D. (2009). Bivariate generalized
exponential distribution. Journal of multivariate analysis, 100(4), 581-593.

\bibitem{36} Lawless, J. F. (1982). Statistical models and methods for
lifetime data. John Wiley and Sons: New York.

\bibitem{30} Lee, H., and Cha, J. H. (2015). On two general classes of
discrete bivariate distributions. The American statistician, 69(3), 221-230.

\bibitem{17} Marsall, A. W., and Olkin, I. (1967). A multivariate
exponential distribution . J. Amer. Statist. Association, 62, 30-44.

\bibitem{26} Mohamed, I., Eliwa, M. S., and El-Morshedy, M., (2017).
Bivariate exponentiated generalized linear exponential distribution with
applications in reliability analysis. https://arxiv.org/abs/1710.00502.

\bibitem{2} Mudholkar, G. S., and Srivastava, D. K. (1993). Exponentiated
Weibull family for analyzing bathtub failure-rate data. IEEE transactions on
reliability, 42, 299--302.

\bibitem{8} Nakagawa, T., and Osaki, S. (1975).The discrete Weibull
distribution. IEEE transactions on reliability, 24(5), 300-301.

\bibitem{14} Nekoukhou, V., Alamatsaz, M.H., and Bidram, H. (2013). Discrete
generalized exponential distribution of a second type. Statistics, 47,
876-887.

\bibitem{15} Nekoukhou,V., and Bidram, H. (2015). The exponentiated discrete
Weibull distribution. SORT, 39(1), 127-146.

\bibitem{31} Nekoukhou, V., and Kundu, D. (2017). Bivariate discrete
generalized exponential distribution. Statistics, 51(5), 1143-1158.

\bibitem{24} Rasool, R., and Akbar, A. J., (2016). On bivariate
exponentiated extended Weibull family of distributions. Ci\^{e}nciae natura,
santa maria, 38(2), 564-576.

\bibitem{10} Roy, D. (2003). The discrete normal distribution.
Communications in statistics - theory and methods, 32(10), 1871-1883.
doi:10.1081/sta-120023256.

\bibitem{11} Roy, D. (2004). Discrete Rayleigh distribution. IEEE
transactions on reliability, 53, 255--260.

\bibitem{20} Sarhan, A., Hamilton, D. C., Smith, B., and Kundu, D., (2011).
The bivariate generalized linear failure rate distribution and its
multivariate extension. Computational statistics and data analysis, 55(1),
644-654.

\bibitem{3} Sarhan, A. M.,and Apaloo, J. (2013). Exponentiated modified
Weibull extension distribution. Reliability engineering and system safety,
112, 137--144.

\bibitem{9} Stein, W. E., and Dattero, R. (1984). A new discrete Weibull
distribution. IEEE transactions on reliability, 33(2), 196-197.

\bibitem{22} Wagner, B. S., and Artur, J. L. (2013). Bivariate Kumaraswamy
distribution: Properties and a new method to generate bivariate classes.
Statistics, 47(6), 1321--1342.

\bibitem{1} Weibull, W. A. (1951). Statistical distribution function of wide
applicability. Journal of applied mechanics, 18, 293--6.

\bibitem{37} Xiao, Jiang., Jeffrey, Chu., and Saralees Nadarajah. (2017).
New classes of discrete bivariate distributions with application to football
data. Communications in statistics - theory and methods, 46:16, 8069-8085,
DOI: 10.1080/03610926.2016.1171358.
\end{thebibliography}
\end{document}